\documentclass[11pt,twoside]{article}

\usepackage{epsfig,amsfonts,color}
\usepackage{amsmath}

\usepackage{amssymb, palatino, geometry,url}
\usepackage{algorithmic}
\usepackage[noresetcount,lined,boxed]{algorithm2e} 

\usepackage[colorlinks=true,linkcolor=blue,citecolor=blue,urlcolor=blue]{hyperref}
\usepackage{subcaption}
\usepackage{amsthm}

\usepackage{fancyhdr}
\newcommand{\btheta}{\boldsymbol{\theta}}

\newcommand{\amin}{\mathop{\mbox{argmin}}}
\newcommand{\amax}{\mathop{\mbox{argmax}}}

\newcommand{\be}{\begin{equation}}
\newcommand{\ee}{\end{equation}}
\newcommand{\bea}{\begin{eqnarray}}
\newcommand{\eea}{\end{eqnarray}}

\newcommand{\bvec}{\left(\begin{array}{c}}
	\newcommand{\evec}{\end{array}\right)}
\newcommand{\bsub}{\begin{subequations}}
	\newcommand{\esub}{\end{subequations}}

\usepackage{lineno}
\usepackage{float}
\usepackage[section]{placeins}
\usepackage[table]{xcolor}

\title{A Computational Framework for Quantifying and \\ Analyzing System Flexibility}
\author{Joshua L. Pulsipher, Daniel Rios, and Victor M. Zavala\thanks{Corresponding Author: victor.zavala@wisc.edu}\\
	{\small Department of Chemical and Biological Engineering}\\
	{\small \;University of Wisconsin, 1415 Engineering Dr, Madison, WI 53706, USA}}
\date{}

\begin{document}
	
\maketitle

\begin{abstract}
We present a computational framework for analyzing and quantifying system flexibility. Our framework incorporates new features that include: general uncertainty characterizations that are constructed using composition of sets,  procedures for computing well-centered nominal points, and a procedure for identifying and ranking flexibility-limiting constraints and critical parameter values. These capabilities allow us to analyze the flexibility of complex systems such as distribution networks. 
\end{abstract}


\section{Introduction and Basic Terminology}
This report details proposed extensions in characterizing the flexibility of physical systems subjected to uncertainties. Such uncertainty may stem from a wide variety of sources including the system's environment (e.g., ambient temperature, product demand) and the system itself (e.g., kinetic constants, battery life). Flexibility denotes the capability of a system to maintain feasible operation over a range of uncertain/random conditions \cite{swaney1985index1}. This becomes prevalent in many applications such as chemical processes \cite{pistikopoulos1995novel,smith2005chemical}, power systems \cite{papaefthymiou2009using,ulbig2015analyzing}, and autonomous vehicles \cite{jun2003path,wang2008autonomous} where it is critical that sufficient flexibility is engineered into system design to promote profitability and safety. For instance, chemical plants encounter numerous uncertainties (e.g., feed flowrates, heat transfer coefficients, equipment fatigue) that can induce costly failures for designs that are not sufficiently flexible \cite{crowl2001chemical}. Accounting for flexibility can also help reduce design costs, since accurate analysis of a system's flexibility helps to prevent the over-engineering that has prevailed traditional chemical process design practices \cite{swaney1985index1}.

A number of approaches have been previously proposed to quantify system flexibility and an extensive review is provided in \cite{grossmann2014evolution}. The stochastic flexibility index, originally proposed by Straub and Grossmann, is a flexibility metric for systems that are subjected to uncertain parameters (e.g., disturbances, physical parameters) that are modeled as random variables \cite{straub1990integrated}. Pistikopoulos and Mazzuchi propose a similar metric in \cite{pistikopoulos1990novel}. The stochastic flexibility index $SF$ is defined as the probability of finding recourse to maintain feasible operation (i.e., the probability of satisfying all the system constraints). This can be computed by integrating the probability density function of the random parameters over the feasible region:

\begin{equation}
	SF := \int_{\btheta\in\Theta} p(\btheta) d\btheta
	\label{eq:SF_definition}
\end{equation}

\noindent where $\boldsymbol{\theta} \in \mathbb{R}^{n_{\theta}}$ is a realization of the random parameters, $\Theta:=\{\btheta\,:\,\psi(\btheta) \leq 0\}$ is the feasible set of the system and $p:\mathbb{R}^{n_{\theta}}\to \mathbb{R}$ is the probability density of the parameters. The function $\psi(\btheta)$ evaluates the feasibility of the system for a given realization of $\btheta$ and is given by:

\begin{equation}
	\begin{aligned}
		& \psi(\btheta) := & \min_{\mathbf{z}, \mathbf{x}} & \max_{j \in J} \ f_j(\mathbf{z}, \mathbf{x}, \btheta) \\
		&& \text{s.t.} & \ \ h_i(\mathbf{z}, \mathbf{x}, \btheta) = 0, & i \in I,
	\end{aligned}
	\label{eq:feasibility_function}
\end{equation}

\noindent where $\mathbf{z} \in \mathbb{R}^{n_z}$ are the system recourse variables, $\mathbf{x} \in \mathbb{R}^{n_x}$ are the system state variables, $f_j(\mathbf{z}, \mathbf{x}, \btheta), \ j \in J$ are the system inequality constraint functions and $h_i(\mathbf{z}, \mathbf{x}, \btheta), \ i \in I$ are the system equality constraint functions. The system is deemed feasible at a given realization of $\btheta$ if $\psi(\btheta) \leq 0$ (this implies that $f_j(\mathbf{z}, \mathbf{x}, \btheta)\leq 0$ and $h_i(\mathbf{z}, \mathbf{x}, \btheta)=0$ hold for all $i\in I$ and $j\in J$).  The system is infeasible at $\btheta$ if $\psi(\btheta)>0$. The boundary of the feasible region is given by the set $\partial \Theta:=\{\btheta\,|\,\psi(\btheta)=0\}$. In other words, at a given realization $\btheta$, the system is at the boundary of the feasible region if $\psi(\btheta)=0$. Problem \eqref{eq:SF_definition} can also be expressed as:
\begin{equation}
	SF = \mathbb{P}\left( \psi(\btheta) \leq 0 \right) = \mathbb{P}\left( \exists \ \mathbf{z}, \mathbf{x} : \max_{j \in J} f_j(\mathbf{z}, \mathbf{x}, \btheta) \leq 0, \  h_i(\mathbf{z}, \mathbf{x}, \btheta) = 0 ,\; i \in I \right).
	\label{eq:SF_definition2}
\end{equation}
From this representation, we can see that the $SF$ index is a joint chance constraint \cite{pulsipher2018MICP}. 

The $SF$ index can be computed rigorously by evaluating  \eqref{eq:SF_definition} via Monte Carlo (MC) sampling. This is done by assessing the feasibility of each realization. Such an approach converges exponentially with the number of MC samples but typically requires a very large number of samples \cite{shapiro2013sample, robert2013monte}. Straub and Grossmann presented a quadrature technique to reduce the number of samples required, but quadrature methods still suffer from scalability issues \cite{straub1990integrated, gerstner1998numerical}. 

\begin{figure}[ht]
	\centering
	\includegraphics[width=0.8\textwidth]{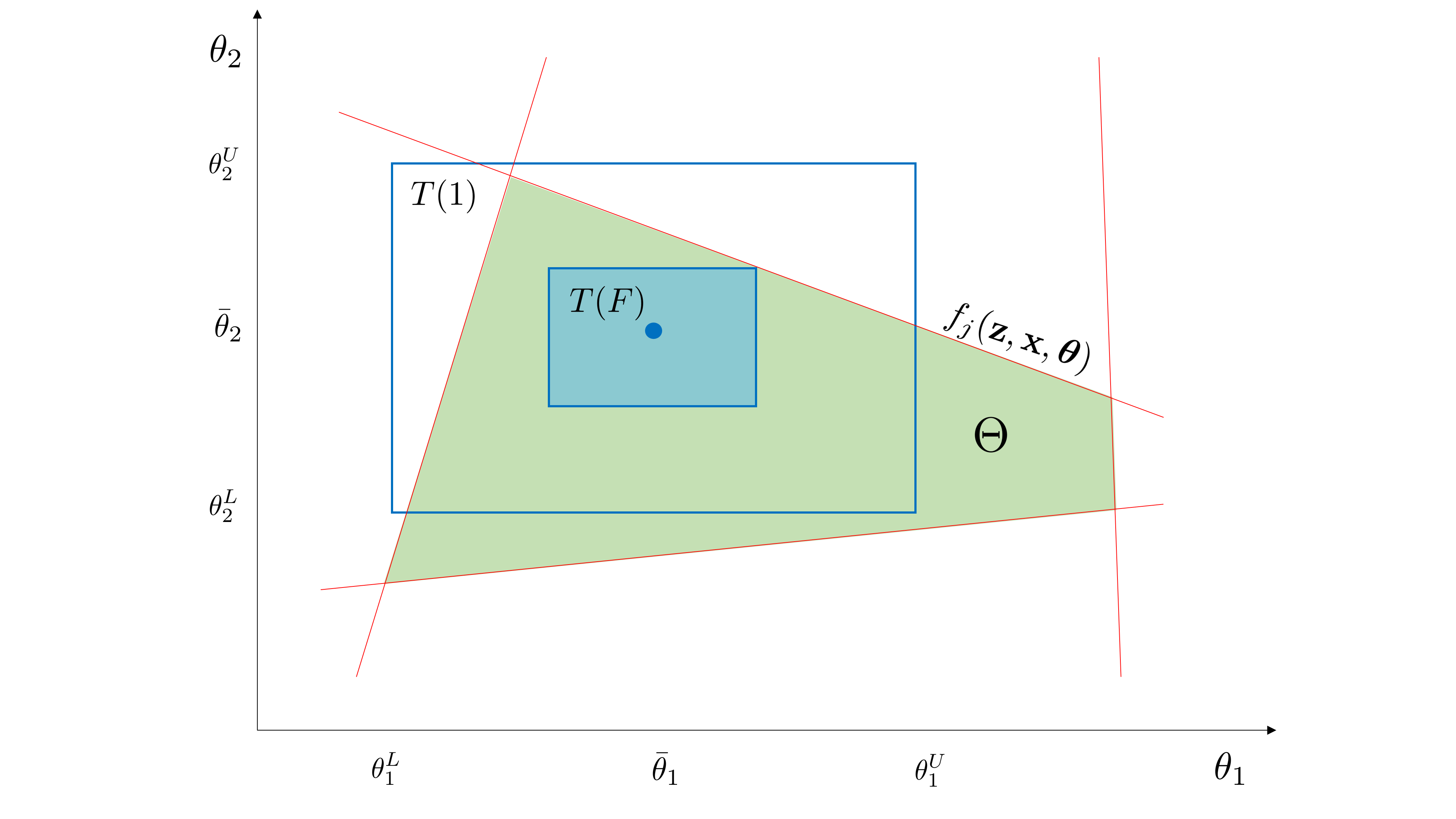}
	\caption{Illustration of flexibility index problem under a hyperbox uncertainty set.}
	\label{fig:flexibility_index}
\end{figure}

The so-called flexibility index problem, first proposed by Grossmann et. al. \cite{grossmann1983optimization}, provides a more scalable approach to quantifying flexibility. This approach seeks to identify the largest uncertainty set $T(\delta)$ (where $\delta \in \mathbb{R}_+$ is a parameter that scales $T$) for which the system remains feasible. In other words, we seek to find the largest uncertainty set under which there exists recourse to recover feasibility. The flexibility index $F$ is defined as:
\begin{equation}
	\begin{aligned}
		& F := & \max_{\delta\in\mathbb{R}_+} &&&\delta \\
		&&\text{s.t.} &&& \max_{\btheta \in T(\delta)} \psi(\btheta) \leq 0.
	\end{aligned}
	\label{eq:flexibility_index_definition}
\end{equation}
This approach is deterministic in nature; consequently, it does not have a direct probabilistic interpretation (as the $SF$ index does). For convenience, $\psi(\btheta)$ can be reformulated by using an upper bounding variable $u \in \mathbb{R}$:
\begin{equation}
	\begin{aligned}
		&\psi(\btheta) = &\min_{\mathbf{z}, \mathbf{x}, u} &&& u \\
		&&\text{s.t.} &&& f_j(\mathbf{z}, \mathbf{x}, \btheta) \leq u, & j \in J \\
		&&&&& h_i(\mathbf{z}, \mathbf{x}, \btheta) = 0, & i \in I.
	\end{aligned}
	\label{eq:feasibility_function2}
\end{equation} 
\noindent Problem \eqref{eq:flexibility_index_definition} is a bi-level optimization problem that seeks to find the largest uncertainty set $T(\delta)$ that can be inscribed in the feasible region defined by $\Theta$. Figure \ref{fig:flexibility_index} illustrates the index $F$ for an inequality system using a hyperbox uncertainty set which is described further below in Equation \eqref{eq:hypercube}.

Swaney and Grossmann \cite{swaney1985index1} proved that Problem \eqref{eq:flexibility_index_definition} is equivalent to searching for the minimum $\delta$ along the boundary $\partial \Theta$, provided that $T(\delta)$ is compact and the constraints $f_j(\mathbf{z}, \mathbf{x}, \btheta)$ and $h_i(\mathbf{z}, \mathbf{x}, \btheta)$ are Lipschitz continuous in $\mathbf{z}$, $\mathbf{x}$, and $\btheta$. Thus, \eqref{eq:flexibility_index_definition} can be expressed as: 

\begin{equation}
	\begin{aligned}
		& F = & \min_{\delta \in \mathbb{R}_+, \ \btheta \in T(\delta)} &&&\delta \\
		&&\text{s.t.} &&& \psi(\btheta) = 0.
	\end{aligned}
	\label{eq:flexibility_index_mod}
\end{equation}

Grossmann and Floudas leverage this observation to propose a mixed-integer program (MIP) that casts the constraint $\psi(\btheta) = 0$ in terms of its first-order Karush-Kuhn-Tucker (KKT) conditions  \cite{grossmann1987active}. The MIP is given by:

\begin{equation}
	\begin{aligned}
		&F = &\min_{\delta, \mathbf{z}, \mathbf{x}, \btheta, \lambda_j, s_j, y_j, \mu_j} &&& \delta \\
		&&\text{s.t.} &&& f_j(\mathbf{z}, \mathbf{x}, \btheta)+ s_j = 0 && j \in J \\
		&&&&& h_i(\mathbf{z}, \mathbf{x}, \btheta) = 0 && i \in I \\
		&&&&& \sum_{j \in J} \lambda_j = 1 \\
		&&&&& \sum_{j \in J} \lambda_j \frac{\partial f_j(\mathbf{z}, \mathbf{x}, \btheta)}{\partial \mathbf{z}} + \sum_{i \in I} \mu_i \frac{\partial h_i(\mathbf{z}, \mathbf{x}, \btheta)}{\partial \mathbf{z}}= 0 \\
		&&&&& s_j \leq U(1 - y_j) && j \in J \\
		&&&&& \lambda_j \leq y_j && j \in J \\
		&&&&& \btheta \in T(\delta)  \\
		&&&&& \lambda_j, s_j \geq 0 ; \ \ \ y_j \in \{0, 1\} && j \in J.
	\end{aligned}
	\label{eq:flex_mip}
\end{equation}

\noindent where $s_j$, $\lambda_j$ are slack variables and Lagrange multipliers for the inequality constraints $f_j(\mathbf{z}, \mathbf{x}, \btheta)$, and $y_j$ are binary variables that indicate if the corresponding inequality constraints are active or inactive. Symbol $\mu_i$ denotes Lagrange multipliers for the equality constraints $h_i(\mathbf{z}, \mathbf{x}, \btheta)$. The constant $U$ is a suitable upper bound for the slack variables. Problem \eqref{eq:flex_mip} provides a tractable formulation to compute the index $F$ that provides a global solution when the convexity conditions presented by Grossmann and Floudas are satisfied (i.e., $\Theta$ is a convex set) \cite{grossmann1987active}. Global solutions for nonlinear systems that do not satisfy these convexity conditions can still be obtained if a global optimization algorithm is applied such as the $\alpha$BB branch-and-bound algorithm employed by Floudas et. al. in \cite{floudas2001global}. The flexibility index problem also shares some common ground with robust optimization methods, Zhang et. al. discuss this relationship in detail for linear systems \cite{zhang2016relation}. 

Because of the difficulty in parameterizing the uncertainty set in terms of a scalar variable $\delta$, the majority of studies reported in the literature use a hyperbox representation of the uncertainty set:
\begin{equation}
	T_{box}(\delta) = \{\btheta : \bar{\btheta} - \delta\Delta\btheta^- \leq \btheta \leq \bar{\btheta} + \delta\Delta\btheta^+\}
	\label{eq:hypercube}
\end{equation}
\noindent where $\bar{\btheta}$ is a {\em nominal point} and $\Delta\btheta^-, \Delta\btheta^+$ are maximum lower and upper deviations, respectively \cite{grossmann2014evolution}. In conjunction with Problem \eqref{eq:flex_mip}, this set provides a straightforward interpretation of the flexibility index (which we denote as $F_{box}$ where the subscript indicates what uncertainty set is used in combination with Problem \eqref{eq:flex_mip}). Specifically, $F_{box} \geq 1$ indicates that the design is sufficiently flexible for the specified deviations. The limitations with the hyperbox set are that it requires upper and lower bounds to be assigned to each random parameters (which may be arbitrary in  applications) and that it does not capture correlations between random parameters. 

Recently, Pulsipher and Zavala \cite{pulsipher2018MICP} showed that any compact set representation of $T(\delta)$ can be utilized in combination with Problem \eqref{eq:flex_mip}. In particular, they showed that an ellipsoidal uncertainty set is the natural choice for $T(\delta)$ when the uncertain parameters are multivariate Gaussian variables $\btheta ~ \sim \mathcal{N}(\bar{\btheta}, V_{\btheta})$ (where $\bar{\btheta}$ is the mean and $V_{\btheta}$ is a positive definite covariance matrix). The ellipsoidal uncertainty set is given by:

\begin{equation}
	T_{ellip}(\delta) = \{\btheta : (\btheta - \bar{\btheta})^T V_{\btheta}^{-1} (\btheta - \bar{\btheta}) \leq \delta \}.
	\label{eq:hyperellipsoid}
\end{equation} 

One advantage of using $T_{ellip}(\delta)$ in combination with Problem \eqref{eq:flex_mip} is that the associated flexibility index (denoted by $F_{ellip}$) can be used to obtain the confidence level:
\begin{equation}
	\alpha^* := \frac{\gamma(\frac{n_{\theta}}{2}, \frac{F_{ellip}}{2})}{\Gamma(\frac{n_{\theta}}{2})}
	\label{eq:alpha_definition}
\end{equation}
\noindent where $\gamma(\cdot)$ and $\Gamma(\cdot)$ are the incomplete and complete gamma functions. Interestingly, this confidence level provides a lower bound for the stochastic flexibility index $SF$ (i.e., $\alpha^*\leq SF)$ \cite{pulsipher2018MICP}. This thus provides an avenue to obtain a probabilistic interpretation of the deterministic flexibility index while avoiding MC sampling and quadrature schemes. 

This paper describes extensions to the aforementioned approaches that help establish a framework for analyzing  flexibility of complex systems. These extensions include a generalization of the uncertainty set representation by using intersecting sets, systematic approaches for obtaining nominal points, and a methodology for identifying and ranking flexibility-limiting constraints. The paper is structured as follows: in Section \ref{sec:extensions} we establish and describe our proposed extensions. Section \ref{sec:examples} illustrates the concepts by using small and large distribution networks. 

\section{Flexibility Analysis Framework} \label{sec:extensions}

This section describes the components of a flexibility analysis framework. The framework incorporates extensions to enable distinct characterizations of $T(\delta)$, to select $\bar{\btheta}$, to compare system designs, and to identify and rank limiting constraints.

\subsection{Uncertainty Set Characterizations} \label{sec:uncert_set}
Problem \eqref{eq:flex_mip} can use any compact set $T(\delta)$ whose size can be parameterized in terms of a scalar variable $\delta$. The selection of $T(\delta)$ will depend on the application and on the nature of the uncertain parameters $\btheta$. In \cite{pulsipher2018MICP}, a number of uncertainty sets are provided that can be used in conjunction with Problem \eqref{eq:flex_mip}. Table \ref{tab:bounded_norms} summarizes such sets. We note that the ellipsoidal norm set in Table \ref{tab:bounded_norms} is equivalent to $T_{ellip}(\delta)$ as defined in \eqref{eq:hyperellipsoid} provided that $A = V_{\btheta}^{-1}$. The ellipsoidal norm set can also be used to approximate polytopes \cite{boyd2004convex}. We also note that $T_{box}(\delta)$ is equivalent to $T_{\infty}(\delta)$ when $\Delta \btheta_i^- = \Delta \btheta_i^+ = 1$. Li et.al. provide an analysis of the $\ell_1$, $\ell_2$, and $\ell_\infty$-norm uncertainty sets along with their mathematical properties in the context of robust optimization \cite{li2011comparative}. The conditional value at risk (CVaR) norm, originally proposed by Pavlikov and Uryasev in \cite{pavlikov2014cvar}, has gained interest in the robust and stochastic optimization communities because of its ability to approximate $\ell_p$ norms via linear programming (precisely reproducing the $\ell_1$ and $\ell_\infty$ norms as extreme cases) \cite{gotoh2016two}. 

\begin{table}[!htb]
	\caption{Potential representations for the uncertainty set $T(\delta)$.}
	\begin{center}
		\begin{tabular}{|c|c|}
			\hline
			Name  & Uncertainty Set  \\ \hline \hline
			Ellipsoidal Norm & $T_{ellip}(\delta) = \{\btheta : ||\btheta - \bar{\btheta}||_{A}^2 \leq \delta \}$ \\
			Hyperbox Set & $T_{box}(\delta) = \{\btheta : \bar{\btheta} - \delta\Delta\btheta^- \leq \btheta \leq \bar{\btheta} + \delta\Delta\btheta^+\}$ \\
			$\ell_\infty$ Norm & $T_{\infty}(\delta) = \{\btheta : ||\btheta - \bar{\btheta}||_\infty \leq \delta \}$ \\
			$\ell_1$ Norm & $T_{1}(\delta) = \{\btheta : ||\btheta - \bar{\btheta}||_1 \leq \delta \}$ \\
			$\ell_2$ Norm & $T_{2}(\delta) = \{\btheta : ||\btheta - \bar{\btheta}||_2 \leq \delta \}$ \\
			CVaR norm & $T_{CVaR}(\delta) = \{\btheta : \langle\langle\btheta - \bar{\btheta}\rangle\rangle_\alpha \leq \delta \}$ \\ \hline
		\end{tabular}
	\end{center}
	\label{tab:bounded_norms}
\end{table}

Li et. al. also analyzed sets that result from intersections (compositions) of $\ell_p$-norm sets (e.g., $T_2(\delta) \cap T_\infty(\delta)$) \cite{li2011comparative}. Inspired by this, we observe that $T(\delta)$ can be represented by any combination of the uncertainty sets in Table \ref{tab:bounded_norms} or other sets (provided that the resulting combined set is compact). This provides the ability to incorporate {\em physical knowledge} on the uncertain parameters and/or to create a wider range of uncertainty set {\em shapes}. For instance, the positive (truncated) ellipsoidal set $T_{ellip+}(\delta) := T_{ellip}(\delta) \cap \mathbb{R}_+^{n_{\theta}}$, where $\mathbb{R}_+^{n_{\theta}} := \{\btheta : \btheta \geq 0\}$ can be used to represent parameters that are known to be non-negative (e.g., demands, prices).  
Similarly, we can define the positive hyperbox set $T_{box+}(\delta) := T_{box}(\delta) \cap \mathbb{R}_+^{n_{\theta}}$. Compositions of sets can be naturally accommodated in conjunction with problem \eqref{eq:flex_mip}, by imposing the constraints corresponding to each set. For instance, $T_{ellip+}(\delta)$ is represented by adding the constraints  $(\btheta - \bar{\btheta})^T V_{\btheta}^{-1} (\btheta - \bar{\btheta}) \leq \delta $ and $\btheta\geq 0$ to the MIP formulation. 

\subsection{Nominal Point Selection Methods}

Most uncertainty sets require a nominal point $\bar{\btheta}$. The choice of the nominal point $\bar{\btheta}$ is critical; for example,  a nominal point that is close to the boundary of the feasible region will have limited flexibility. In some applications, specifying $\bar{\btheta}$ might be difficult if not enough data is available or if a system is not routinely operated at any given point (e.g., a power grid). Thus, we introduce methods that can be employed to find well-centered nominal points.  We note that one would be tempted to let $\bar{\btheta}$ be a free variable in Problem \eqref{eq:flex_mip} to find a nominal point but this approach leads to trivial solutions. This is because the objective is minimized when the nominal point $\bar{\btheta}$ is placed on the boundary of $\btheta$.  Consequently, this naive approach cannot be used to compute a nominal point. 

We propose an extension of the {\em analytic center} (commonly used in interior-point methods) to compute a nominal point. The analytic center $\bar{\btheta}_{ac}$ is given by:
\begin{equation}
	\begin{aligned}
		&{\bar{\btheta}}_{ac} := & \amax_{\btheta, \mathbf{z}, \mathbf{x}, s_j} &&& \sum_{j \in J} \log \left(s_j \right) \\
		&&\text{s.t.} &&& f_j(\mathbf{z}, \mathbf{x}, \btheta) +s_j\leq 0, && j \in J \\
		&&&&& h_i(\mathbf{z}, \mathbf{x}, \btheta) = 0, && i \in I.
	\end{aligned}
	\label{eq:analytical_center}
\end{equation}
Here, $s_j\in \mathbb{R}$ are slack variables. This problem pushes the constraints to the interior of the feasible region. Any bounded solution of the above problem satisfies $s_j> 0$ and thus $f_j(\mathbf{z},\mathbf{x},\bar{\btheta}_{ac})<0$ for all $j\in J$. Consequently, we have that $\psi(\bar{\btheta}_{ac})<0$ and thus $\bar{\btheta}_{ac}\in \Theta$ (the analytic center is feasible).   We also note that the analytic center is the point that maximizes the geometric mean of the constraints (slack variables). By changing the objective to $\sum_{j\in J}s_j $ one finds a center point that maximizes the arithmetic mean.

We also propose a center point $\bar{\btheta}_{fc}$ (that we call the {\em feasible center}) and that we define as:
\begin{equation}
	\begin{aligned}
		&\bar{\btheta}_{fc} := &\amax_{\btheta, \mathbf{z}, \mathbf{x}, s} &&& s \\
		&&\text{s.t.} &&& f_j(\mathbf{z}, \mathbf{x}, \btheta) + s \leq 0, & j \in J \\
		&&&&& h_i(\mathbf{z}, \mathbf{x}, \btheta) = 0, & i \in I.
	\end{aligned}
	\label{eq:feasible_center}
\end{equation} 
Here, $s\in \mathbb{R}$ is a slack variable. The feasible center is the point that maximizes the worst-case slack variable (as opposed to the geometric mean used in the analytic center). We also note that $\bar{\btheta}_{fc} = \amin_{\btheta}\;\psi(\btheta)=\amax_{\btheta}\;-\psi(\btheta)$ and that $\psi(\bar{\btheta}_{fc})\leq 0$.  Moreover,  \eqref{eq:feasible_center} is equivalent to \eqref{eq:feasibility_function2} when the parameter $\theta$ is set as a free variable. 

The feasible center is an approximation of the center point: 
\begin{align}
\bar{\btheta}_d:= \amax_{\btheta\in \Theta}\; d(\btheta,\partial \Theta) 
\end{align}
where $d(\btheta,\partial\Theta)=\min_{\btheta'\in \partial\Theta}\; \|\btheta-\btheta'\|$. In other words, $\bar{\btheta}_d$ is the point that maximizes the distance to the boundary of feasible set.  This point is a natural  center point but it is difficult to compute due to the nested max and min operators. We now proceed to show that $ |-\psi(\btheta)|$ provides a lower bound of $d(\btheta,\partial \Theta)$.  Consequently, because the feasible center $\bar{\btheta}_{fc}$ maximizes $-\psi(\btheta)$, we have that this also implicitly maximizes $d(\btheta,\partial \Theta)$.  This property can be established as follows, if the constraints are Lipschitz continuous (as often assumed in the flexibility analysis literature), we have that $\psi(\cdot)$ is Lipschitz as well and thus:
\begin{align}
|\psi(\hat{\btheta})-\psi(\btheta)|\leq L\|\bar{\btheta}_d-\btheta\|
\end{align}
where $\hat{\btheta} = \amin_{\btheta'\in \partial\Theta}\; \|\btheta-\btheta'\|$ for a particular value of $\btheta$. Because $\hat{\btheta}\in \partial \Theta$, we have that $\psi(\hat{\btheta})=0$ and thus,
\begin{align}
|\psi(\hat{\btheta})-\psi(\btheta)|&=|-\psi(\btheta)|\nonumber\\
&\leq L\|\hat{\btheta}-\btheta\|\nonumber\\
&=L \min_{\btheta'\in \partial \Theta}\; \|\btheta'-\btheta\|\nonumber\\
&=L \,d(\btheta,\partial\Theta). 
\end{align}
We thus conclude that:
\begin{align}
 d(\btheta,\partial\Theta)\geq \frac{1}{L} |-\psi(\btheta)|.
\end{align}
Consequently, maximizing $-\psi(\btheta)$ implicitly maximizes $d(\btheta,\partial\Theta)$. Moreover, we note that $d(\btheta,\partial\Theta)=-\psi(\btheta)=0$ if $\btheta\in \partial \Theta$.

\subsection{Design Comparison} \label{sec:flex_compare}

The flexibility index is a valuable metric that can be used to compare system designs. The utility and interpretation of the index for a system will depend on the choice of $T(\delta)$. For instance, the index can be used to determine the confidence level $\alpha^*$ if the set $T_{ellip}(\delta)$ is used (when $\btheta ~ \sim \mathcal{N}(\bar{\btheta}, V_{\btheta})$). The index might not have a clear interpretation or utility for some choices of $T(\delta)$ but can still be useful to quantify {\em improvements in flexibility} from design modifications (retrofits). We illustrate this feature by using the following example. 

Consider two system designs that are each subjected to two normal random variables and that have no recourse variables $\bf{z}$ as summarized in Table \ref{tab:2d_constrs}.
\begin{table}[!htb]
	\caption{System Constraints of Design A and Design B.}
	\begin{center}
		\begin{tabular}{|r|c|c|}
			\hline
			 & Design A & Design B  \\ \hline \hline
			 $f_1:$ & $\btheta_1 + \btheta_2 - 14 \leq 0$ & $\frac{3}{4}\btheta_1 + \btheta_2 - 14 \leq 0$ \\
			 $f_2:$ & $\btheta_1 - 2 \btheta_2 - 2 \leq 0$ & $\frac{3}{4}\btheta_1 - 2 \btheta_2 - 2 \leq 0$ \\
			 $f_3:$ & $-\btheta_1 \leq 0$ & $-\btheta_1 \leq 0$ \\
			 $f_4:$ & $-\btheta_2 \leq 0$ & $-\btheta_2 \leq 0$ \\ \hline
		\end{tabular}
	\end{center}
	\label{tab:2d_constrs}
\end{table}

\noindent The nominal point $\bar{\btheta}$ and the covariance matrix $V_{\btheta}$ are taken to be:

\begin{equation}
	\bar{\btheta} =
	\begin{bmatrix}
	4 \\ 5
	\end{bmatrix}
	\ \ \ \ \ \
	V_{\btheta} = 
	\begin{bmatrix}
	2 & 1 \\
	1 & 3
	\end{bmatrix}.
	\label{eq:2d_example_conditions}
\end{equation}

\noindent Given that the parameters are multivariate Gaussian, we choose to use $T_{ellip}(\delta)$. We also compute $F_{box}$ using $T_{box}(\delta)$, where $\Delta\btheta^-$ and $\Delta\btheta^+$ are both taken to be $(4.243, 5.196)$. These bounds correspond to $\bar{\btheta}_i \pm 3\sigma_i$ confidence bounds where $\sigma_i$ is the standard deviation of $\btheta_i$. For comparison, the $SF$ index is rigorously computed using 100,000 MC samples. Table \ref{tab:2d_results} summarizes the results. Both flexibility indexes indicate that Design B is more flexible. We also see that the confidence level $\alpha^*$ also provides a conservative approximation of $SF$. 

\begin{table}[!htb]
	\caption{Comparing the flexibility of Design A and Design B.}
	\begin{center}
		\begin{tabular}{|c|cccc|}
			\hline
			& $F_{box}$ & $F_{ellip}$ & $\alpha^*$ (\%) & $SF$-MC (\%) \\ \hline \hline
			Design A & 0.53 & 3.57 & 83.1 & 96.6\\
			Design B & 0.66 & 6.40 & 95.9 & 98.9 \\ \hline
		\end{tabular}
	\end{center}
	\label{tab:2d_results}
\end{table}

\begin{figure}[!htb]
	\centering
	\begin{subfigure}[t]{.45\textwidth}
		\centering
		\includegraphics[width=\textwidth]{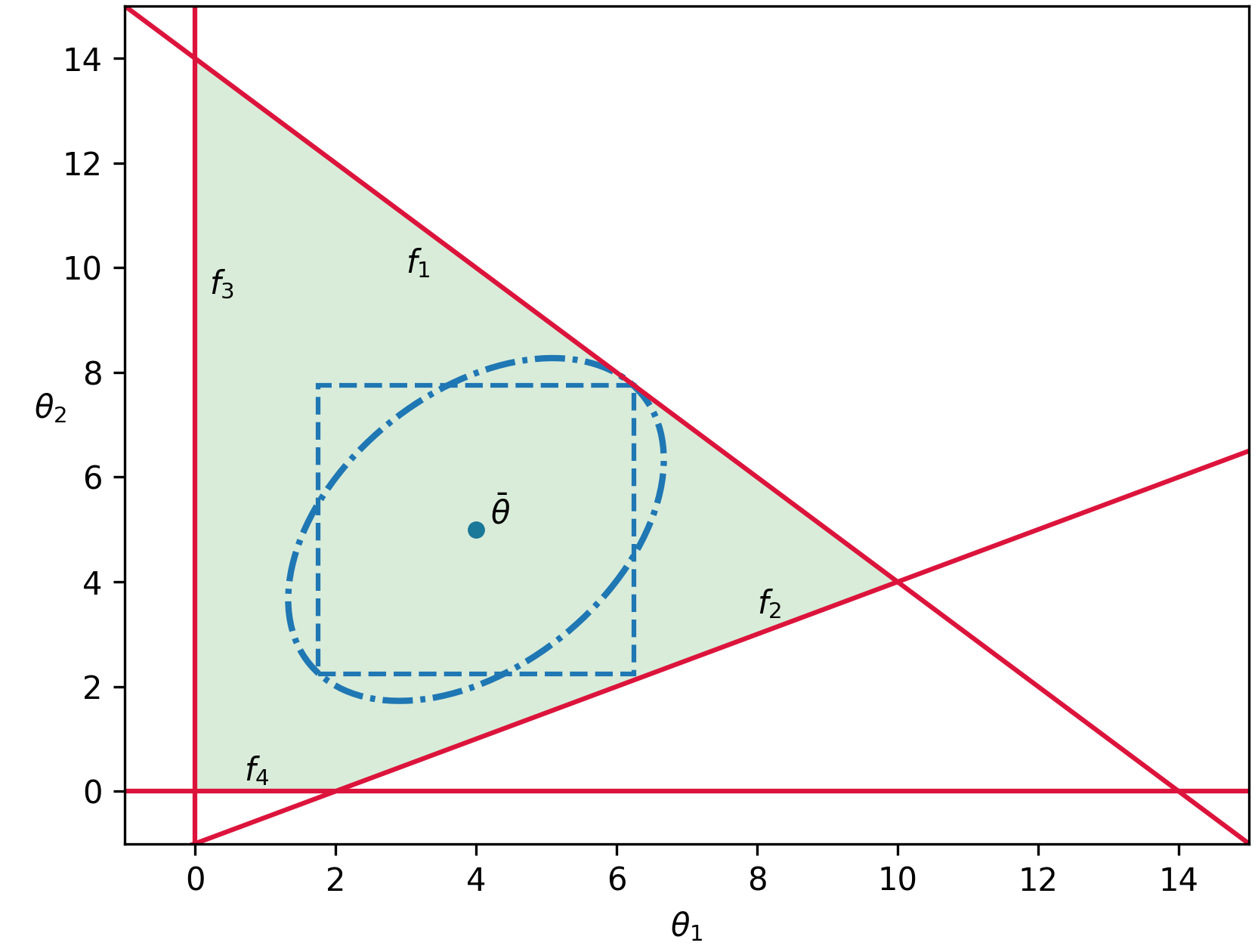}
		\caption{Design A}
	\end{subfigure}
	\quad
	\begin{subfigure}[t]{.45\textwidth}
		\centering
		\includegraphics[width=\textwidth]{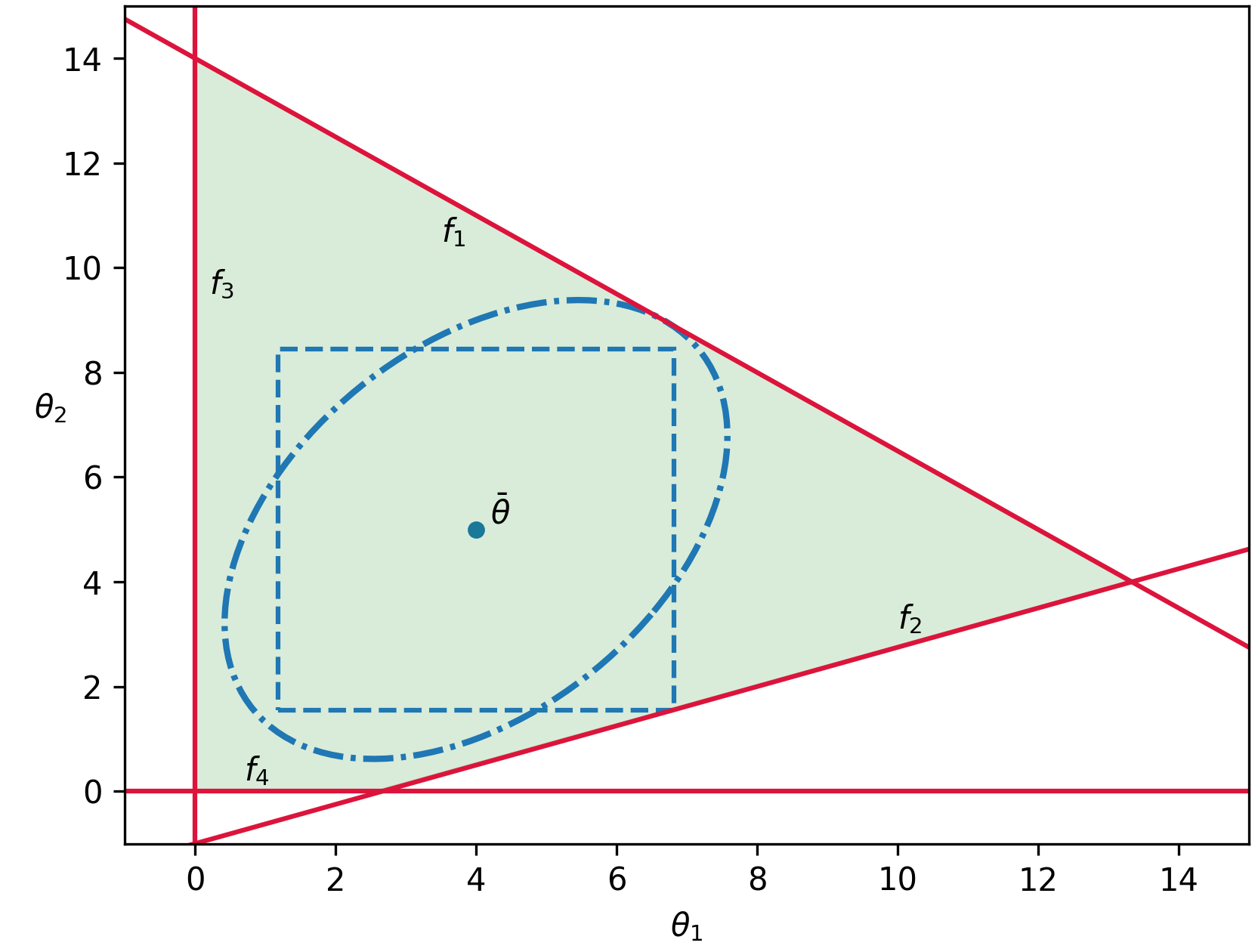}
		\caption{Design B}
	\end{subfigure}
	\caption{Elliptical and hyperbox uncertainty sets relative to feasible regions of Designs A and B.}
	\label{fig:2d_compare_example}
\end{figure}

Figure \ref{fig:2d_compare_example} depicts the two different systems. It is apparent that Design B has a larger feasible region $\Theta$ and thus has a larger $SF$ index. A key observation is that the choice of $T(\delta)$ significantly affects how it measures system flexibility. For instance, for Design B, the ellipsoidal and hyperbox sets identify different limiting constraints, and thus exhibit distinct behavior. In particular, index $F_{box}$ is limited by constraint $f_2$ and would not change if constraints $f_1$, $f_3$, and/or $f_4$ were varied to increase the area of $\Theta$, even though the overall system flexibility would be improved. On the other hand, index $F_{ellip}$ is limited by constraint $f_1$. This conservative behavior parallels what is observed with the use of uncertainty sets in robust optimization, where such sets are used to optimize against the worst case scenario such that the solution is robust in the face of uncertainty, and the conservativeness of a solution depends on the chosen shape of the uncertainty set \cite{gorissen2015practical}. 

\subsection{Limiting Constraint Ranking} \label{sec:flex_rank}

An important observation that we make is that the flexibility index problem also implicitly identifies inequality constraints that limit flexibility. Such information can be valuable in design or operating procedures. Specifically, at a solution of Problem \eqref{eq:flex_mip}, the binary variables $y_j$ indicate which inequality constraints are active and therefore limit flexibility. Problem \eqref{eq:flex_mip} can thus be resolved by excluding this first set of limiting constraints, so that the next set of limiting constraints can be identified. This step can be repeated to rank system constraints to a desired extent (as long as the solution of Problem \eqref{eq:flex_mip} is bounded). This is done by noticing that each set of limiting constraints will have an associated flexibility index, which indicates how limiting a particular set of constraints is relative to the other sets. We note that this analysis operates under the assumption that the inequality constraints can be eliminated/relaxed such that the system still has physical meaning. However, this assumption is met in many applications since constraints that cannot be removed (e.g., physical laws) typically correspond to equality constraints.

This methodology can also be used to identify and rank {\em system components} that limit flexibility in order to guide design improvements or quantify value of specific components, since such identification is useful to understand and identify system vulnerabilities. Thus, this ranking becomes useful for complex systems where it is difficult to determine vulnerable components. For example, this approach can identify which heat-exchangers in a heat-exchanger network could be retrofitted to most increase the system flexibility, and similarly identify heat-exchangers that would not lead to increased flexibility when improved. Limiting constraint information can also be used to quantify the impact of {\em failure} of system components (e.g., a production facility). 

This proposed approach shares some common ground with optimal design problems where design variables $\mathbf{d}$ (added to the system constraints $f_j(\mathbf{z}, \mathbf{x}, \btheta)$ and $h_i(\mathbf{z}, \mathbf{x}, \btheta)$) are selected to minimize cost and satisfy a specified index $F$ \cite{bansal2000flexibility,grossmann1979optimum,pistikopoulos1988optimal}. This provides a straightforward approach for selecting a system design that is sufficiently flexible prior to its determination. However, our approach differs in that it can be used to fully characterize how each component of a particular system design impacts system flexibility instead of choosing design variables to improve flexibility to desired extent. This can becomes useful for intelligently retrofitting existing systems and for understanding which components most limit flexibility (i.e., are most vulnerable to uncertainties) that therefore warrant close monitoring.

We illustrate the constraint ranking methodology using Design A (described in Section \ref{sec:flex_compare}). Here, we choose the set $T_{ellip}(\delta)$. Problem \eqref{eq:flex_mip} is solved to identify the limiting constraints, which are then eliminated to solve the problem again to identify a new set of limiting constraints. In this case, the problem is solved four times as only one constraint becomes active at the solution of each problem. The results are summarized in Table \ref{tab:2d_compare_results}. We observe that constraint $f_1$ limits the flexibility index the most. Specifically, the value of $F_{ellip}$ corresponding to the next limiting constraint $f_2$ is 79.3\% larger. Diminishing changes in $F_{ellip}$ are obtained for subsequent limiting constraints, showing that constraints $f_2$-$f_4$ have little impact on the system flexibility relative to $f_1$. 

\begin{table}[!htb]
	\caption{Constraint ranking results for Design A.}
	\begin{center}
		\begin{tabular}{|c|ccc|}
			\hline
			       & Active Constraint & $F_{ellip}$ & Flexibility Increase (\%) \\ \hline \hline
			Rank 1 & $f_1$             & 3.57        & $-$         \\
			Rank 2 & $f_2$             & 6.40        & 79.3        \\
			Rank 3 & $f_3$             & 8.00        & 124.1       \\
			Rank 4 & $f_4$             & 8.40        & 135.3 \\ \hline
		\end{tabular}
	\end{center}
	\label{tab:2d_compare_results}
\end{table}

\begin{figure}[!htb]
	\centering
	\begin{subfigure}[t]{.45\textwidth}
		\centering
		\includegraphics[width=\textwidth]{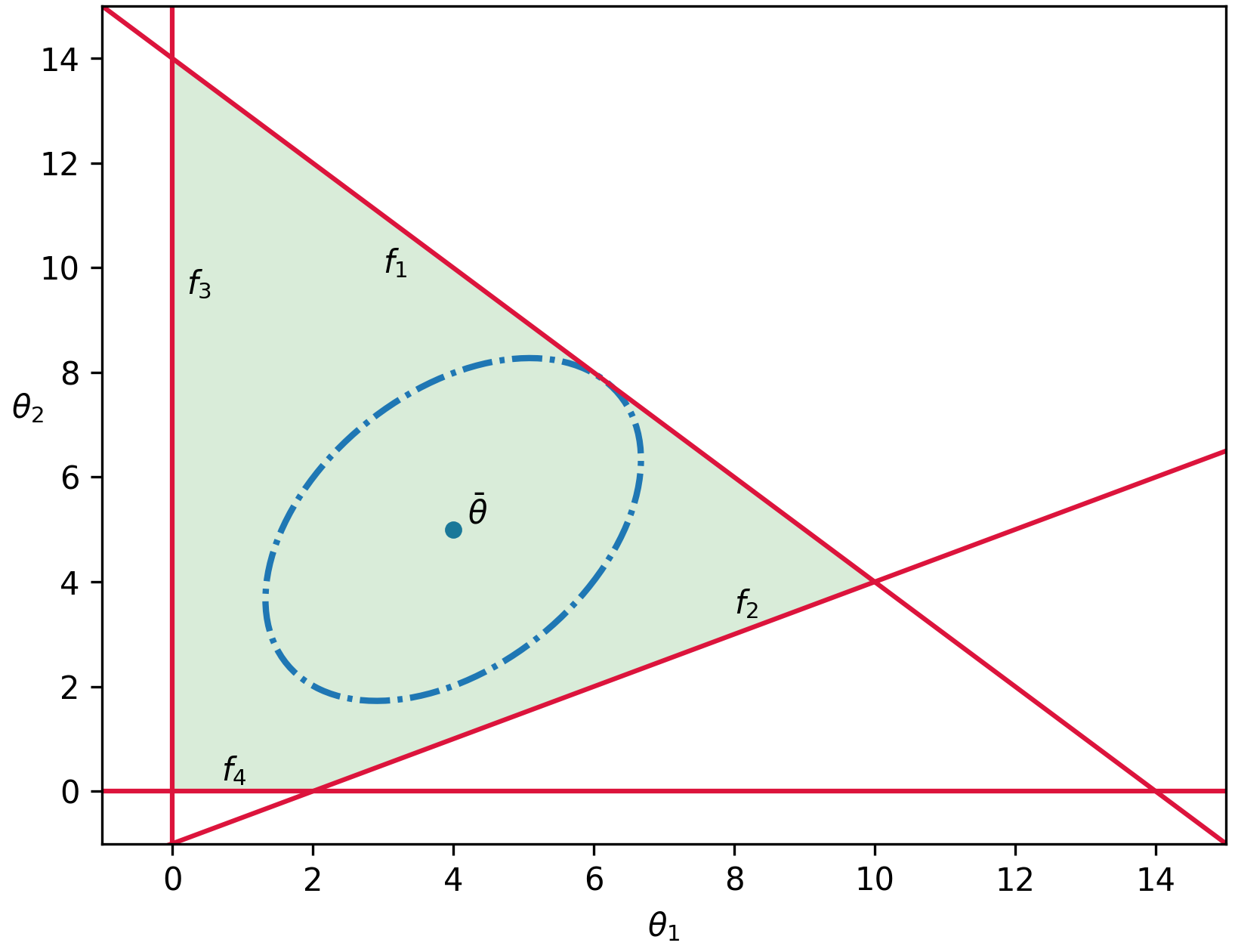}
		\caption{Rank 1: Constraint $f_1$ is limiting}
	\end{subfigure}
	\quad
	\begin{subfigure}[t]{.45\textwidth}
		\centering
		\includegraphics[width=\textwidth]{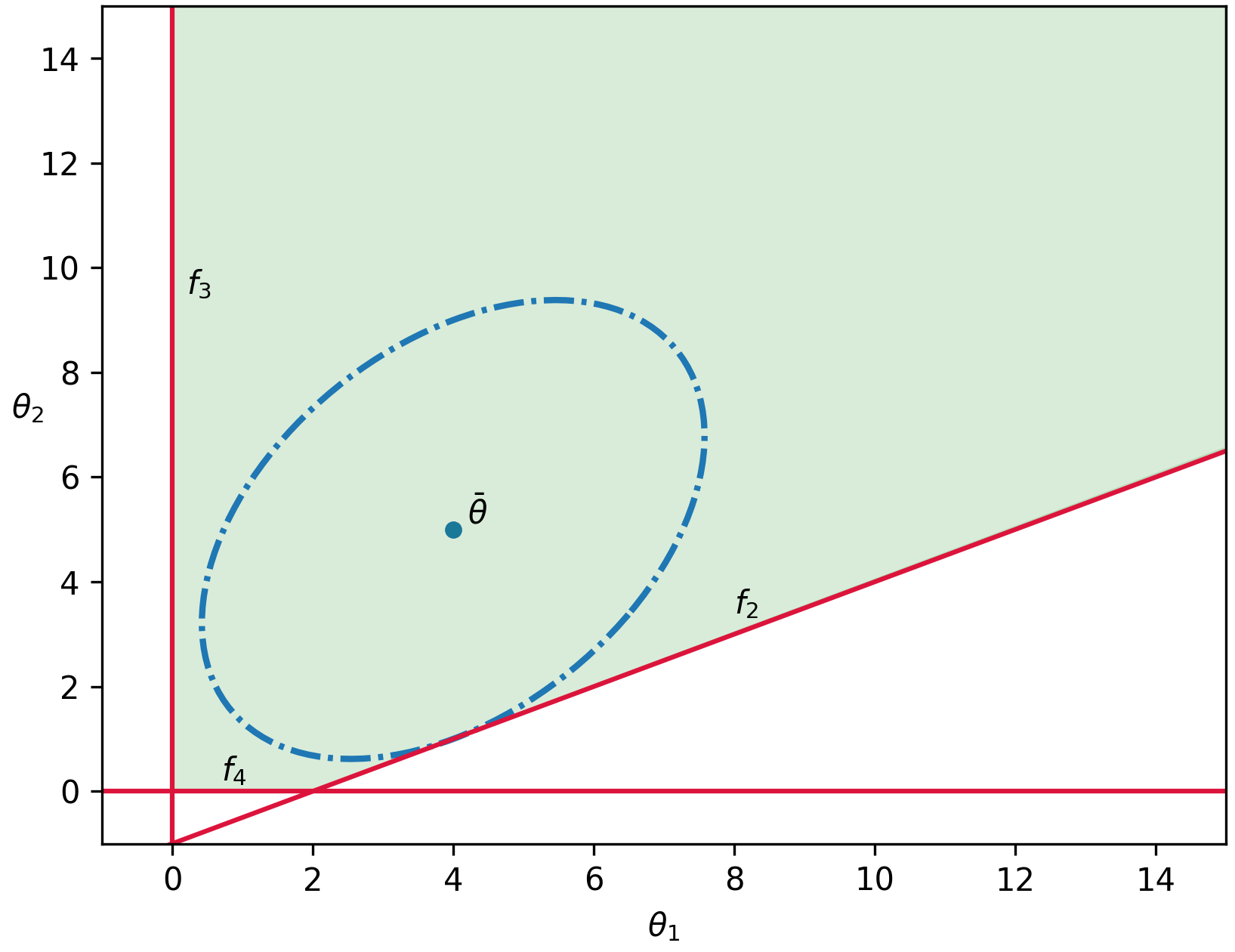}
		\caption{Rank 2: Constraint $f_2$ is limiting}
	\end{subfigure} \\
	\vspace{1em}
	\begin{subfigure}[t]{.45\textwidth}
		\centering
		\includegraphics[width=\textwidth]{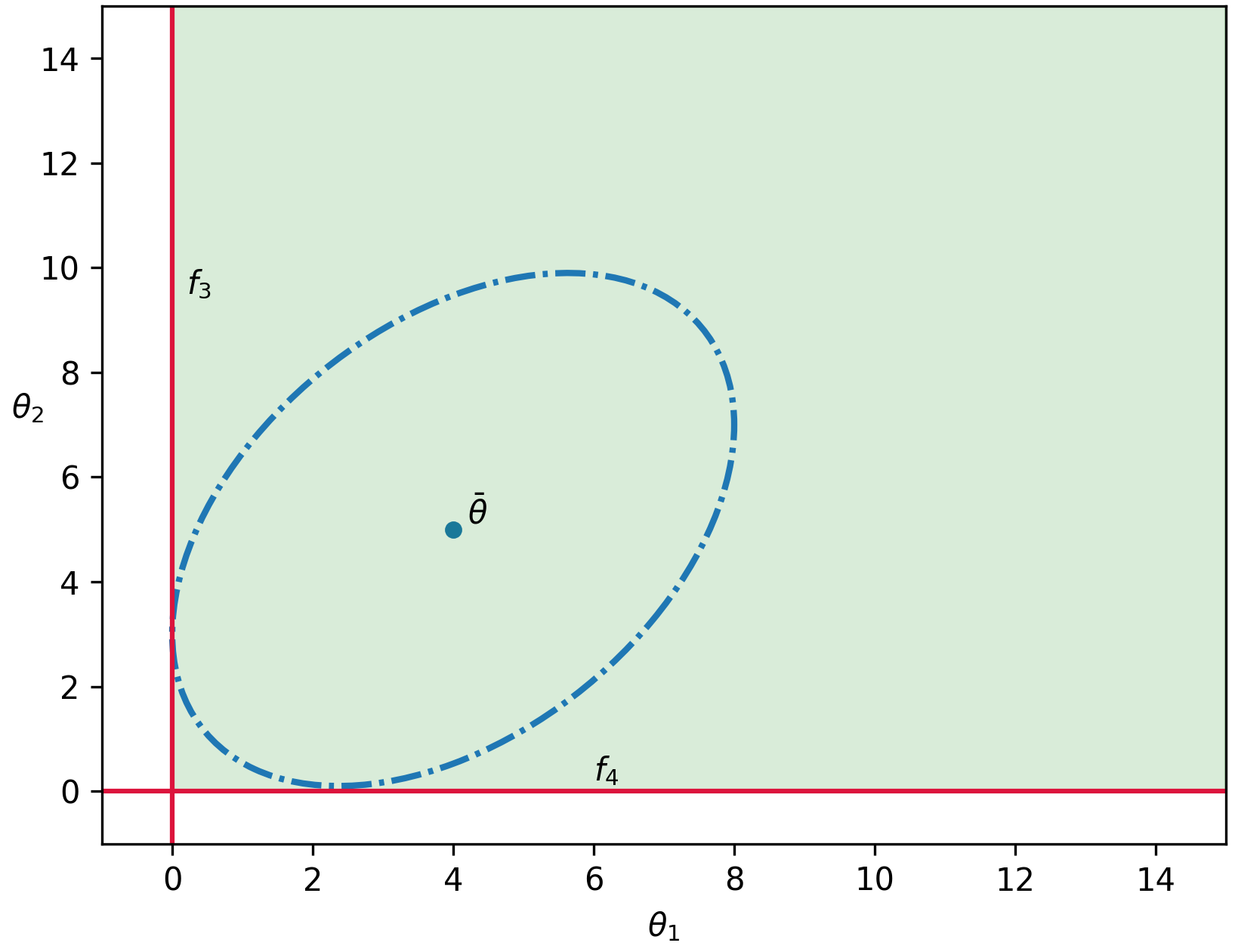}
		\caption{Rank 3: Constraint $f_3$ is limiting}
	\end{subfigure}
	\quad
	\begin{subfigure}[t]{.45\textwidth}
		\centering
		\includegraphics[width=\textwidth]{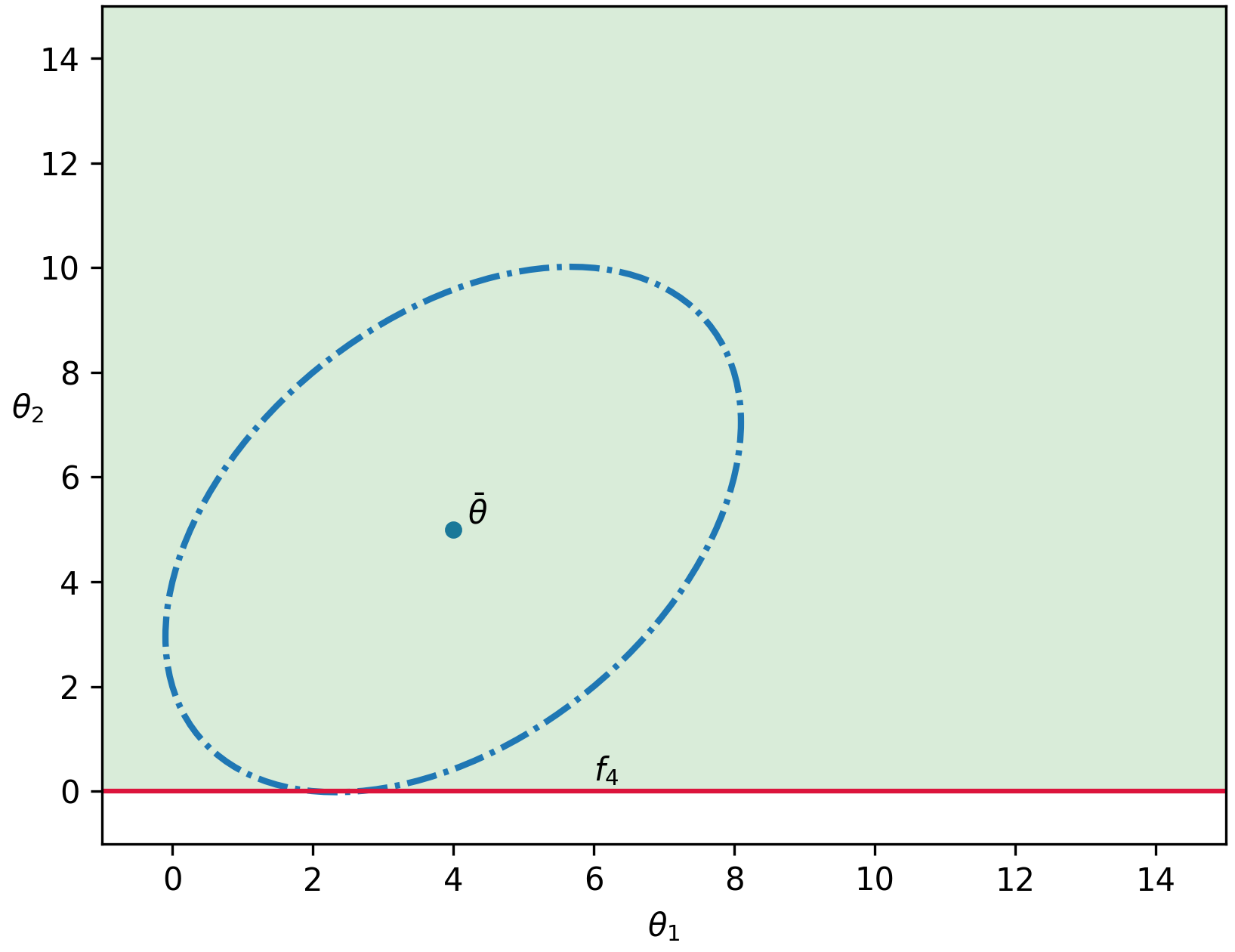}
		\caption{Rank 4: Constraint $f_4$ is limiting}
	\end{subfigure}
	\caption{Identifying and ranking limiting constraints using an elliptical uncertainty set.}
	\label{fig:ranking_example}
\end{figure}

Figure \ref{fig:ranking_example} depicts the solutions corresponding to the results presented in Table \ref{tab:2d_compare_results}. It is apparent that the uncertainty set is larger after each subsequent ranking step. We again make the observation that the limiting constraints are determined by the shape of the uncertainty set. Thus, different uncertainty sets can yield distinct classifications of limiting constraints. 

\FloatBarrier

\section{Case Study: Distribution Networks}\label{sec:examples}

We apply the flexibility analysis framework to distribution networks. We consider a simple three-node network and power distribution networks from standard libraries (IEEE-14 and 141-Bus Network). These results seek to highlight interesting insights that can be gained with the framework. The data used in the power networks is obtained from the MATPOWER package \cite{zimmerman2011matpower}. All of the formulations are implemented in FlexJuMP v0.1.0 (\url{https://github.com/zavalab/FlexJuMP.jl}), a Julia package we have developed that extends JuMP \cite{DunningHuchetteLubin2017} to streamline the implementation of the flexibility analysis framework presented in this paper. These formulations are solved using Pavito to leverage Gurobi 7.5.2 in combination with IPOPT for efficient solution of general convex MIPs or using Pajarito to leverage Gurobi 7.5.2 for efficient solution of conic MIPs on an Intel\textregistered \, Core\texttrademark \, i7-7500U machine running at 2.90 GHz with 4 hardware threads and 16 GB of RAM.

\subsection{Physical Model}
 Each network is modeled by performing balances at each node $n \in \mathcal{C}$ and enforcing capacity constraints on the arcs $a_k, k \in \mathcal{A},$ and on the suppliers $s_b, b \in \mathcal{S}$. The demands $d_m, m \in \mathcal{D},$ are assumed to be the uncertain parameters. The flexibility index problem thus seeks to identify the largest set of simultaneous demand withdrawals that the system can tolerate. The deterministic network model is given by:
\begin{subequations}
	\begin{equation}
		\sum_{k \in \mathcal{A}_n^{rec}} a_k - \sum_{k \in \mathcal{A}_n^{snd}} a_k + \sum_{b \in \mathcal{S}_n} s_b - \sum_{m \in \mathcal{D}_n} d_m = 0, \ \ \ n \in \mathcal{C}
		\label{eq:node_balances}
	\end{equation}
	\begin{equation}
		-a_k^C \leq a_k \leq a_k^C, \ \ \ k \in \mathcal{A}
		\label{eq:arc_constrs}
	\end{equation}
	\begin{equation}
		0 \leq s_b \leq s_b^C, \ \ \ b \in \mathcal{S}
		\label{eq:supp_constrs}
	\end{equation}
	\label{eq:network_model}
\end{subequations}
\noindent where $\mathcal{A}_n^{rec}$ denotes the set of receiving arcs at node $n$, $\mathcal{A}_n^{snd}$ denotes the set of sending arcs at $n$, $\mathcal{S}_n$ denotes the set of suppliers at $n$, $\mathcal{D}_n$ denotes the set of demands at $n$, $a_k^C$ are the arc capacities, and $s_b^C$ are the supplier capacities. The Lagrange multipliers associated with the constraints in \eqref{eq:arc_constrs} are denoted as $\lambda^L_k$ (corresponding to the lower arc bounds) and $\lambda^U_k$ (corresponding to the upper bounds). Similarly, we define multipliers corresponding to the supply constraints in \eqref{eq:supp_constrs} as $\gamma^{L}_b$ and $\gamma^{U}_n$.

\subsection{Three-Node Network} \label{sec:3node}
We consider three-node distribution network designs (see Figure \ref{fig:3_node_designs}). Design 1 corresponds to a centralized supplier base case \cite{zavala2017stochastic}. Design 2 differs from Design 1 in that it employs a decentralized supply scheme, and Design 3 differs from Design 1 in that it contains an additional transportation arc. Each design is subjected to multivariate Gaussian demands $\btheta = (d_1,d_2,d_3) \sim \mathcal{N}(\bar{\btheta}, V_{\btheta})$. The nominal point $\bar{\btheta}$ is set to be the analytic and feasible centers. The center points obtained are given by $\bar{\btheta}_{ac}=(0.0, 50.0, 0.0)$ and $\bar{\btheta}_{fc}=(0.0, 37.4, 0.0)$. The covariance matrix $V_{\btheta}$ is assumed to be:
\begin{equation}
	V_{\btheta} = 
	\begin{bmatrix}
		80 & \beta & \beta \\
		\beta & 80 & \beta \\
		\beta & \beta & 120
	\end{bmatrix}
	\label{eq:3d_covars}
\end{equation}
\noindent where $\beta \in \{-40, 0, 50\}$ is a parameter that captures covariance. These choices of $\bar{\btheta}$ and $V_{\btheta}$ are used in conjunction with four types of uncertainty sets: $T_{ellip}(\delta)$, $T_{ellip+}(\delta)$, $T_{box}(\delta)$, and $T_{box+}(\delta)$. The hyperbox deviations $\Delta\btheta^-, \Delta\btheta^+$ are both taken to be $(26.833, 26.833, 32.863)$, which correspond to $\bar{\btheta}_i \pm 3\sigma_i$ confidence bounds. 
Detailed results for all uncertainty sets and designs are provided in Appendix \ref{a:sup_mat}. The optimality of each solution was verified by evaluating feasibility of 10,000 points generated randomly along the optimized uncertainty set boundaries. The stochastic flexibility index was also estimated for each instance by using 1,000,000 MC samples. 

\begin{figure}[!htb]
	\centering
	\begin{subfigure}[t]{.31\textwidth}
		\centering
		\includegraphics[width=\textwidth]{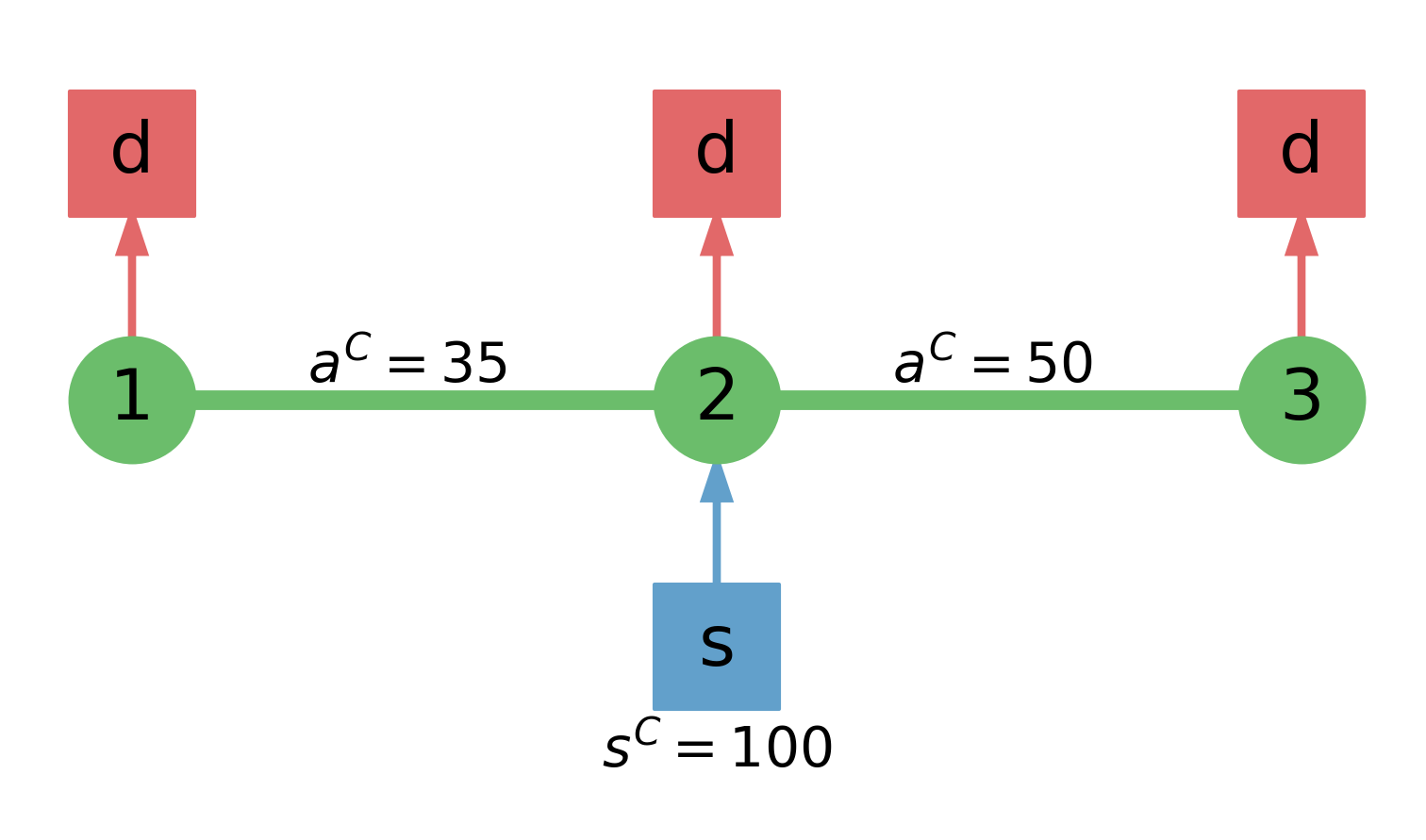}
		\caption{Design 1}
	\end{subfigure}
	\quad
	\begin{subfigure}[t]{.31\textwidth}
		\centering
		\includegraphics[width=\textwidth]{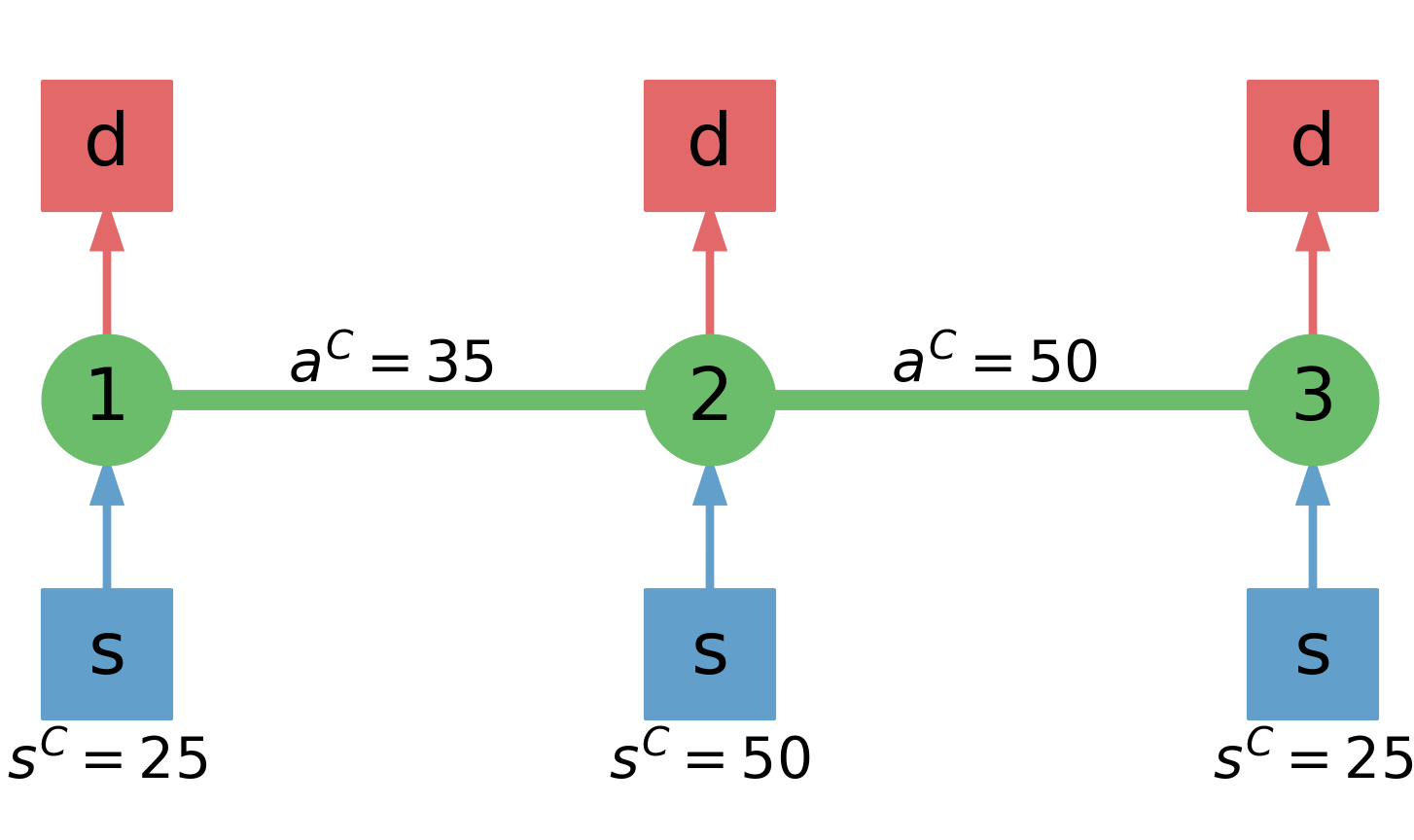}
		\caption{Design 2}
	\end{subfigure}
	\quad
	\begin{subfigure}[t]{.31\textwidth}
		\centering
		\includegraphics[width=\textwidth]{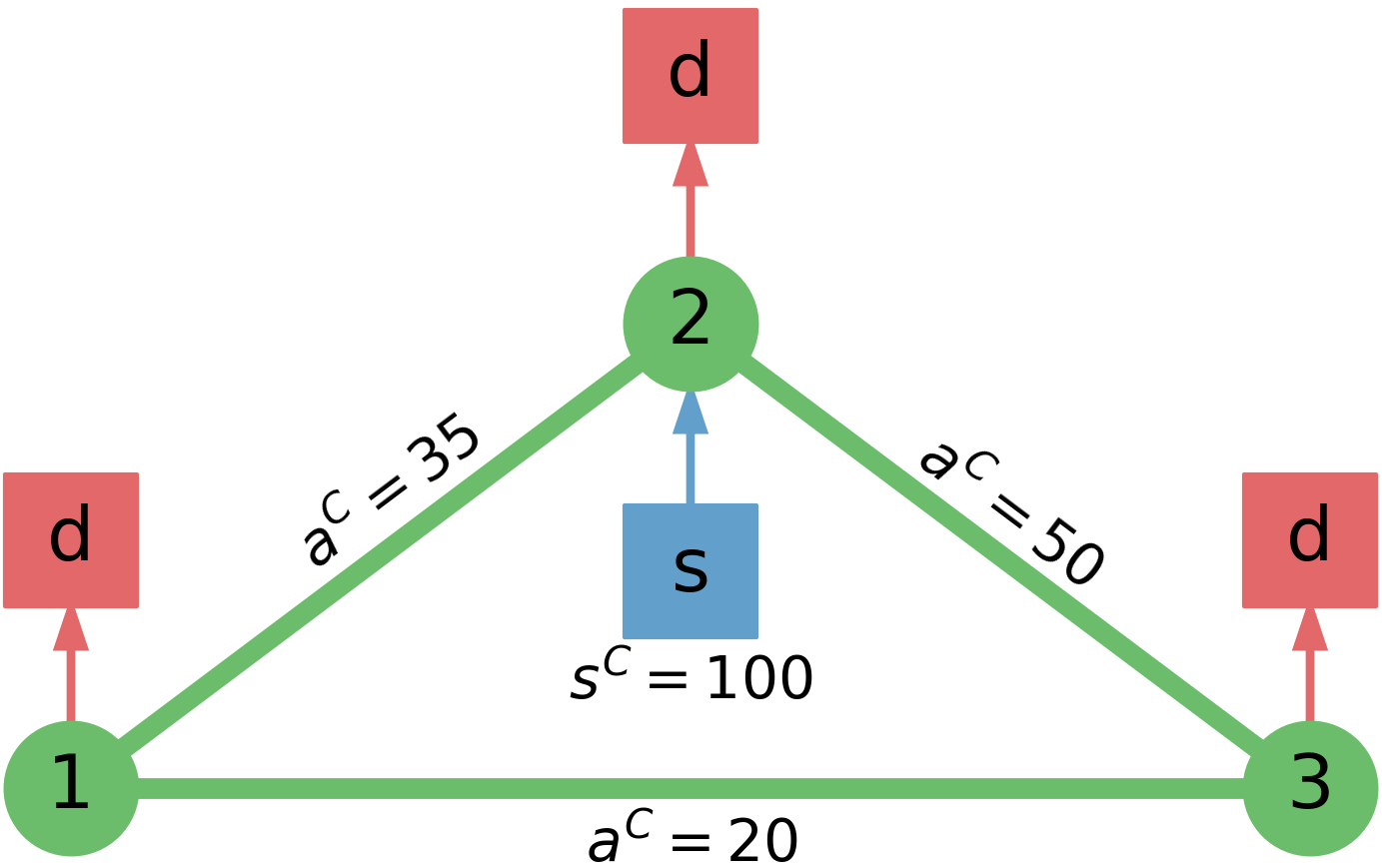}
		\caption{Design 3}
	\end{subfigure}
\end{figure}

\subsubsection{Uncertainty Set and Nominal Point Comparison}

To illustrate the differences between the various uncertainty sets and center points, we consider those used in conjunction with Design 1. Table \ref{tab:3node_results} summarizes these results for the instances with $\beta = 0$. The value of $SF$ corresponding to the feasible center is smaller than that obtained with the analytic center under multivariate Gaussian uncertainty; whereas the converse is true under the truncated multivariate Gaussian. Furthermore, we observe that the behavior of all the indexes tested are in agreement with the trend observed in the stochastic index $SF$. The interpretation of  the flexibility index is fundamentally different for the ellipsoidal and hyperbox sets; however, we observe that both follow the same trends in this instance because they have the same set of limiting constraints. We also note that computing each of these indexes required less than 0.2 seconds in contrast to the MC-generated index $SF$ (which required around 5,200 seconds on average for each case). Also, the hyperbox sets required less computing time than the ellipsoidal sets (this is expected since Problem \eqref{eq:flex_mip} is linear with the hyperbox sets and conic with the ellipsoidal sets).  

\begin{table}[!htb]
	\caption{Flexibility index results for Design 1 of three-node distribution network with $\beta = 0$.}
	\begin{center}
		\begin{tabular}{|c|ccccc|}
			\hline
			Center   & Set                  & $F$   & Active Constraint   & $SF$-MC (\%) & Solution Time (s) \\ \hline \hline
			         & $T_{ellip}(\delta)$  & 8.93  & $\gamma_{2}^L$         & 99.71        & 0.18 \\
			analytic & $T_{ellip+}(\delta)$ & 8.93  & $\gamma_{2}^U$         & 99.44        & 0.13 \\
			         & $T_{box}(\delta)$    & 0.578 & $\gamma_{2}^L$         & 99.71        & 0.03 \\
			         & $T_{box+}(\delta)$   & 0.578 & $\gamma_{2}^U$         & 99.44        & 0.04 \\ \hline
			         & $T_{ellip}(\delta)$  & 5.00  & $\gamma_{2}^L$         & 98.71        & 0.15 \\
			feasible & $T_{ellip+}(\delta)$ & 14.00 & $\gamma_{2}^U$         & 99.95        & 0.09 \\
			         & $T_{box}(\delta)$    & 0.432 & $\gamma_{2}^L$         & 98.71        & 0.03 \\
			         & $T_{box+}(\delta)$   & 0.723 & $\gamma_{2}^U$         & 99.95        & 0.03 \\ \hline
		\end{tabular}
	\end{center}
	\label{tab:3node_results}
\end{table}

\begin{figure}[!htb]
	\centering
	\begin{subfigure}[t]{.48\textwidth}
		\centering
		\includegraphics[width=\textwidth]{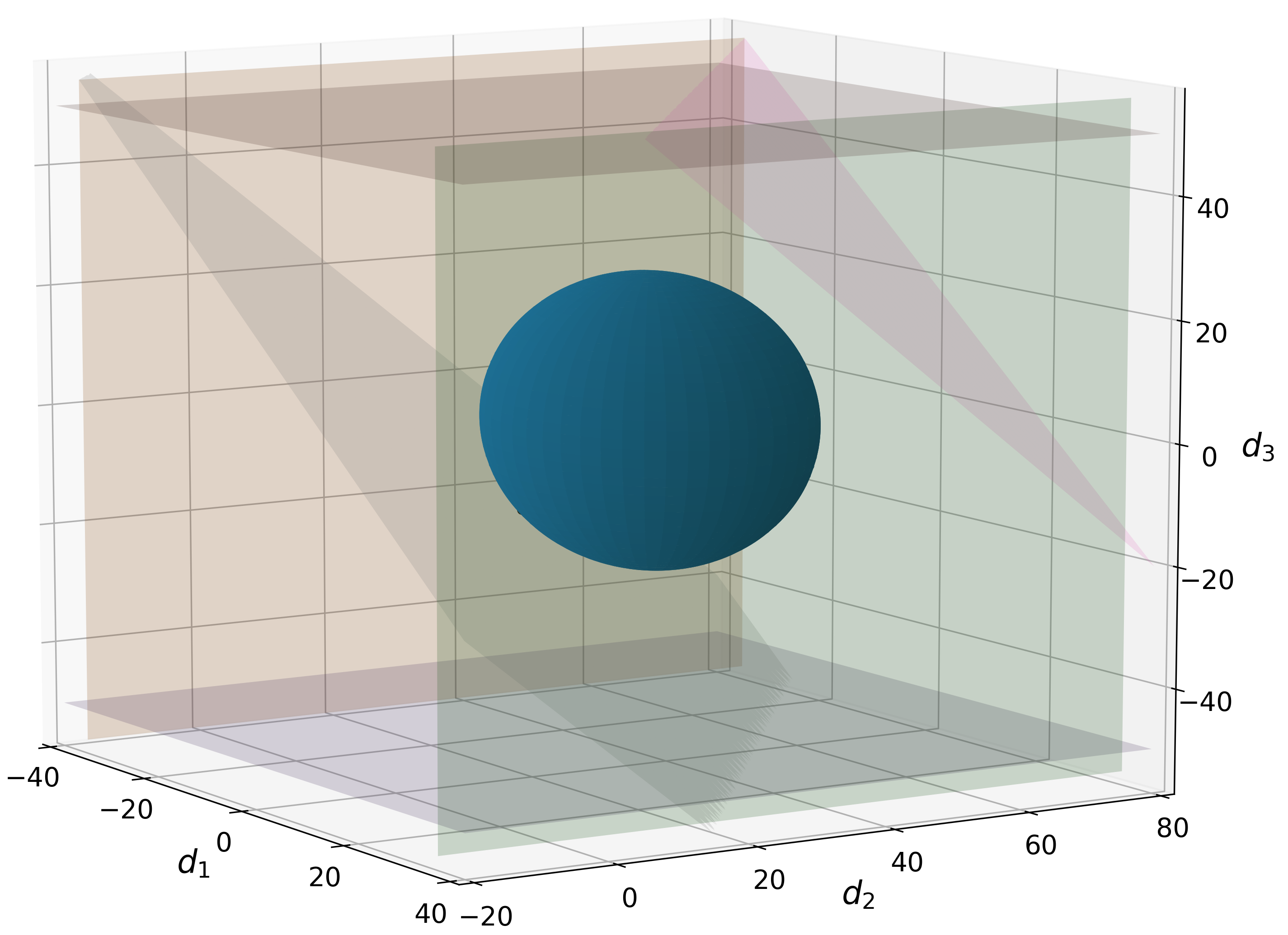}
		\caption{$T_{ellip}(\delta)$}
		\label{fig:ellip_fc}
	\end{subfigure}
	\quad
	\begin{subfigure}[t]{.48\textwidth}
		\centering
		\includegraphics[width=\textwidth]{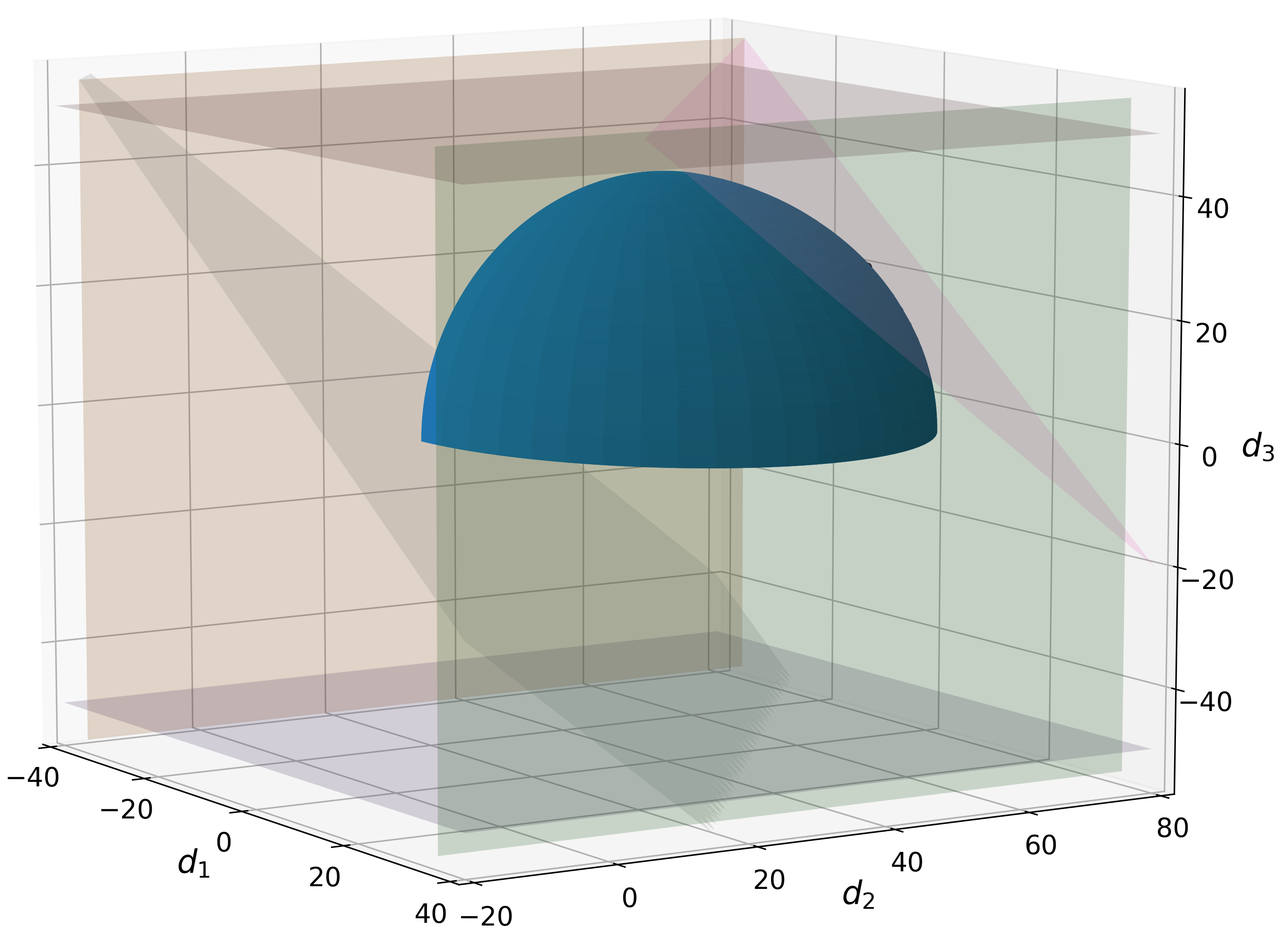}
		\caption{$T_{ellip+}(\delta)$}
	\end{subfigure} \\
	\vspace{1em}
	\begin{subfigure}[t]{.48\textwidth}
		\centering
		\includegraphics[width=\textwidth]{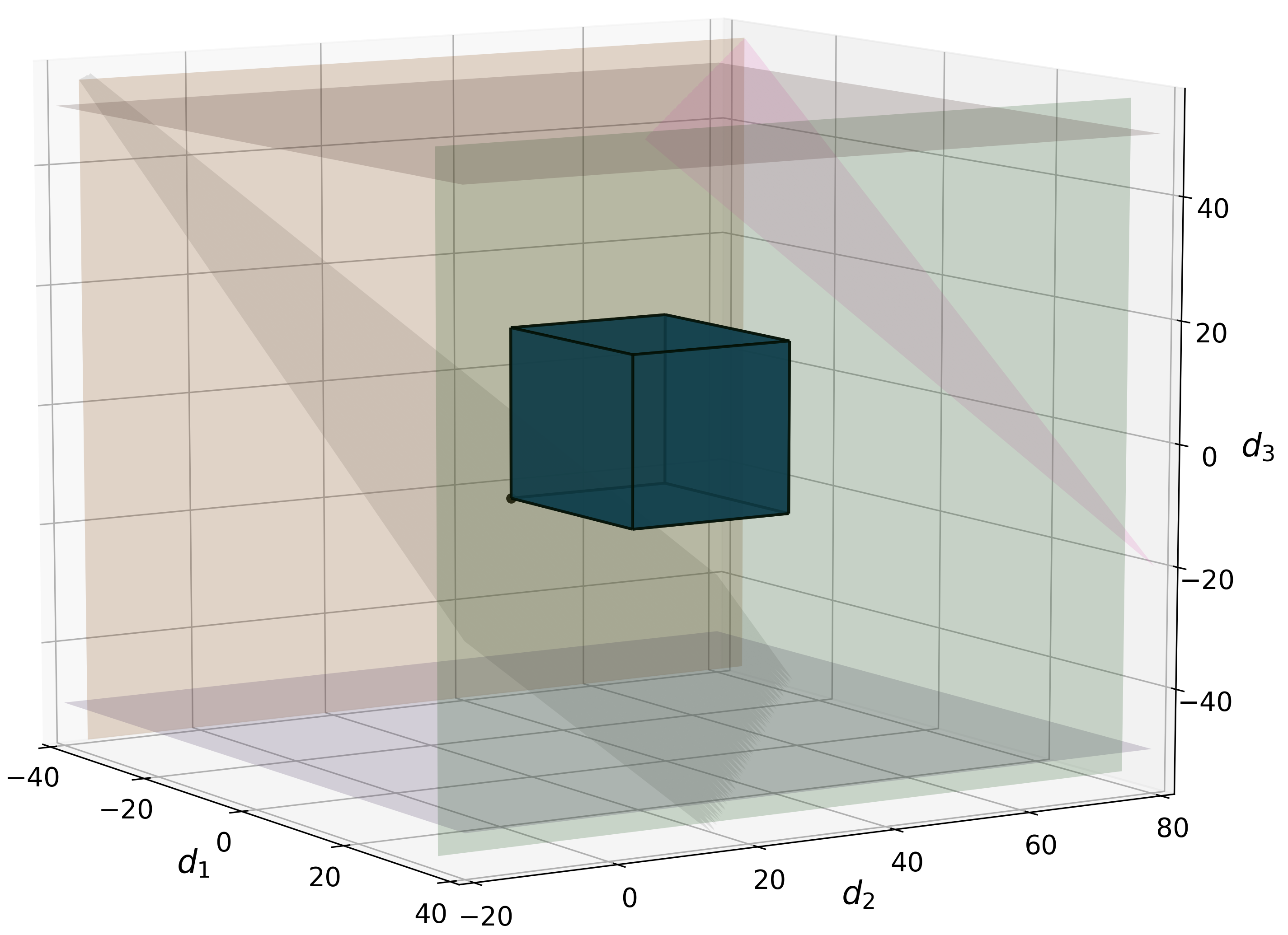}
		\caption{$T_{box}(\delta)$}
	\end{subfigure}
	\quad
	\begin{subfigure}[t]{.48\textwidth}
		\centering
		\includegraphics[width=\textwidth]{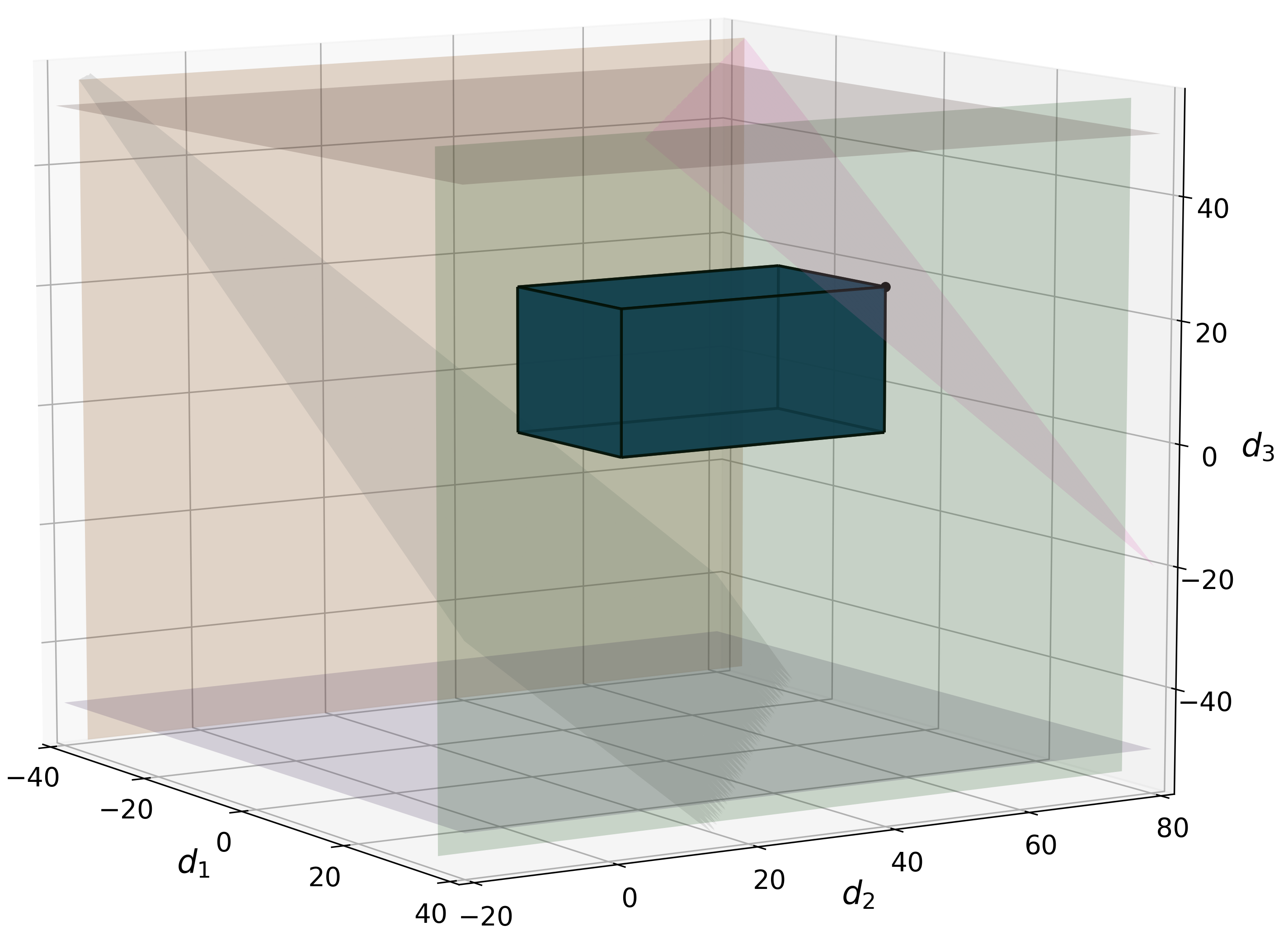}
		\caption{$T_{box+}(\delta)$}
	\end{subfigure}
	\caption{A graphical representation of the optimized ellipsoidal and hyperbox uncertainty sets with $\bar{\btheta} = \bar{\btheta}_{fc}$ corresponding to Design 1 of the 3-node network.}
	\label{fig:3_node_sets}
\end{figure}

\begin{figure}[!htb]
	\centering
	\begin{subfigure}[t]{.48\textwidth}
		\centering
		\includegraphics[width=\textwidth]{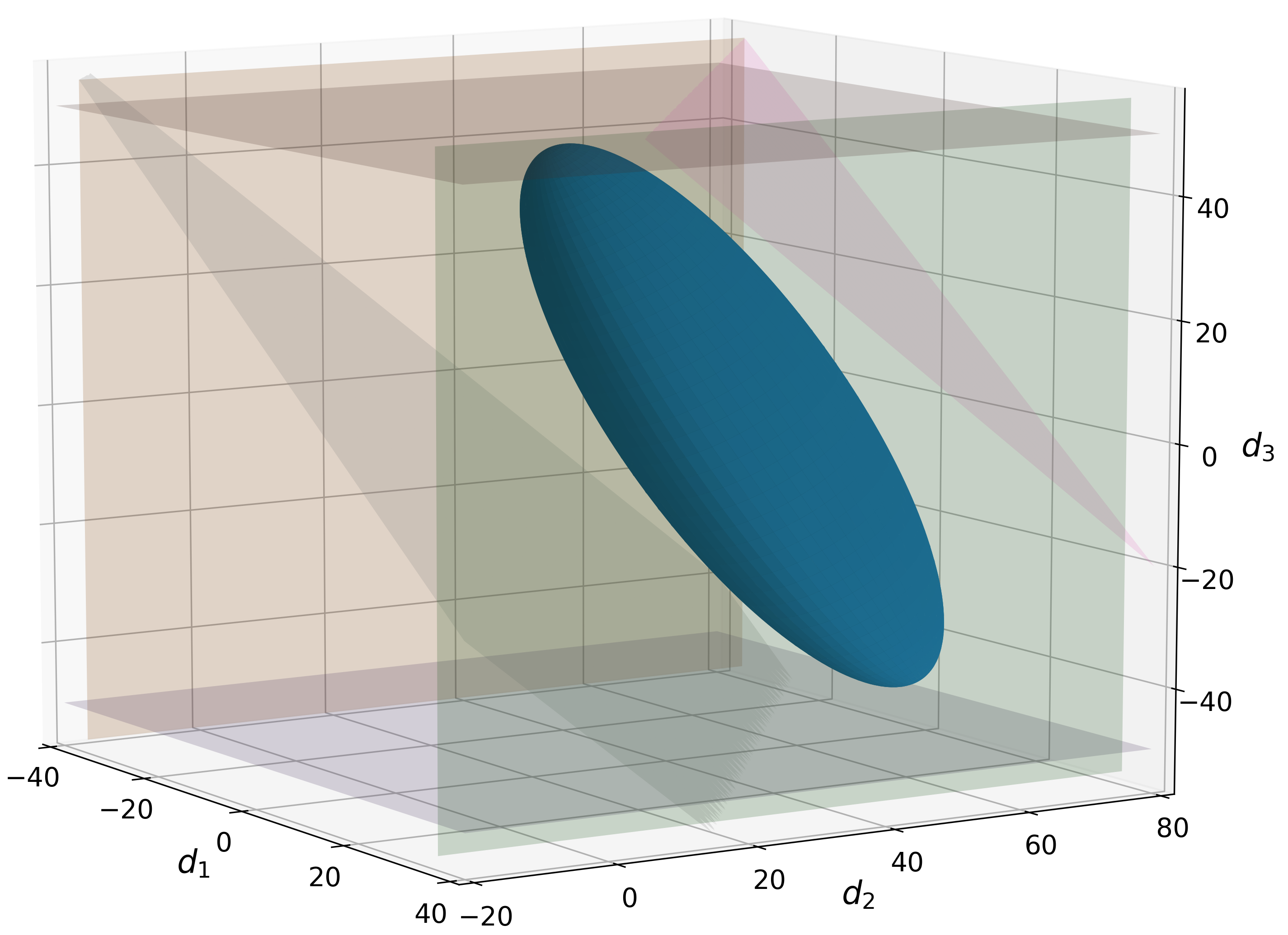}
		\caption{$\beta = -40$}
	\end{subfigure}
	\quad
	\begin{subfigure}[t]{.48\textwidth}
		\centering
		\includegraphics[width=\textwidth]{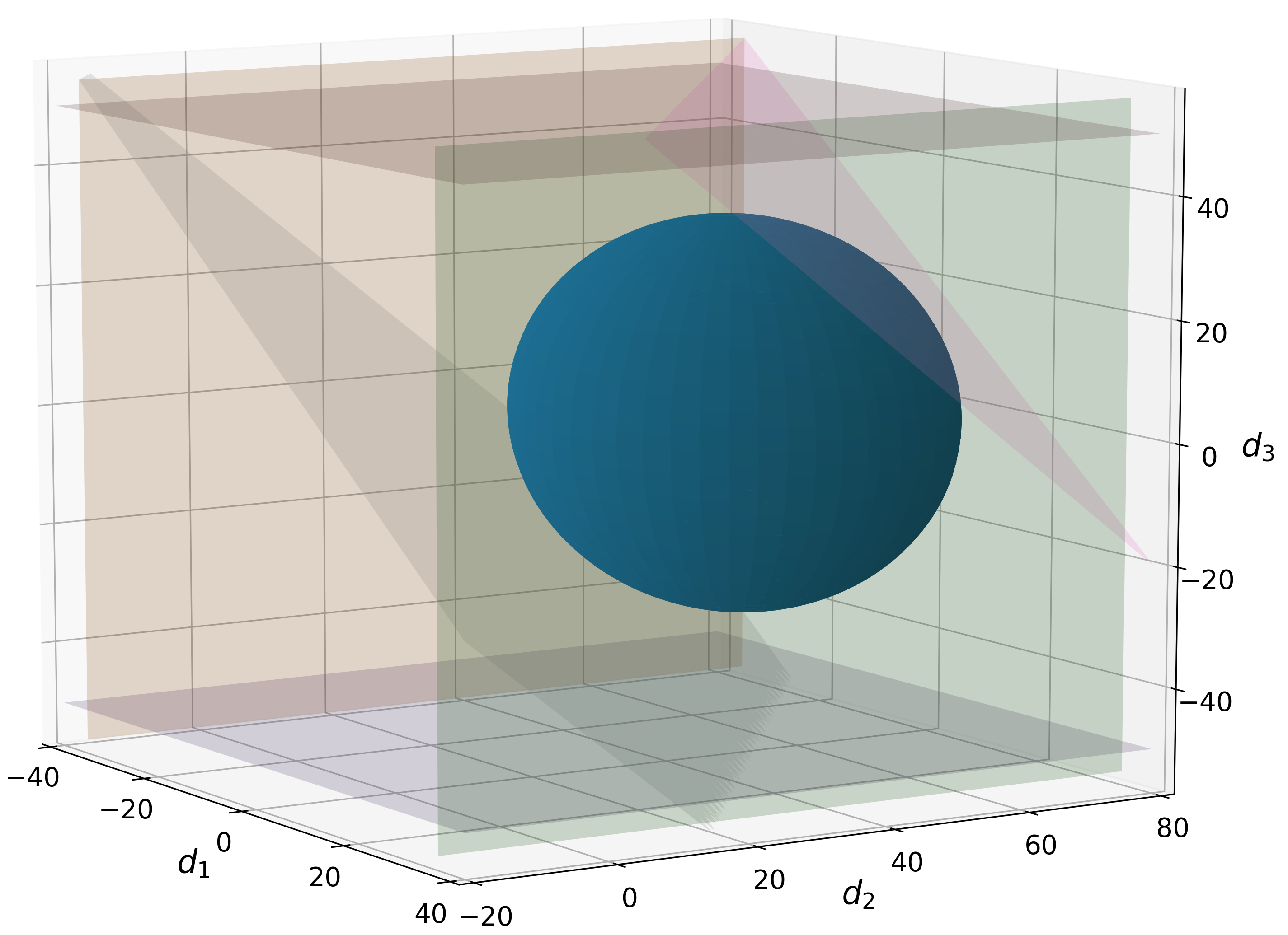}
		\caption{$\beta = 0$}
		\label{fig:ellip_ac}
	\end{subfigure} \\
	\vspace{1em}
	\begin{subfigure}[t]{.48\textwidth}
		\centering
		\includegraphics[width=\textwidth]{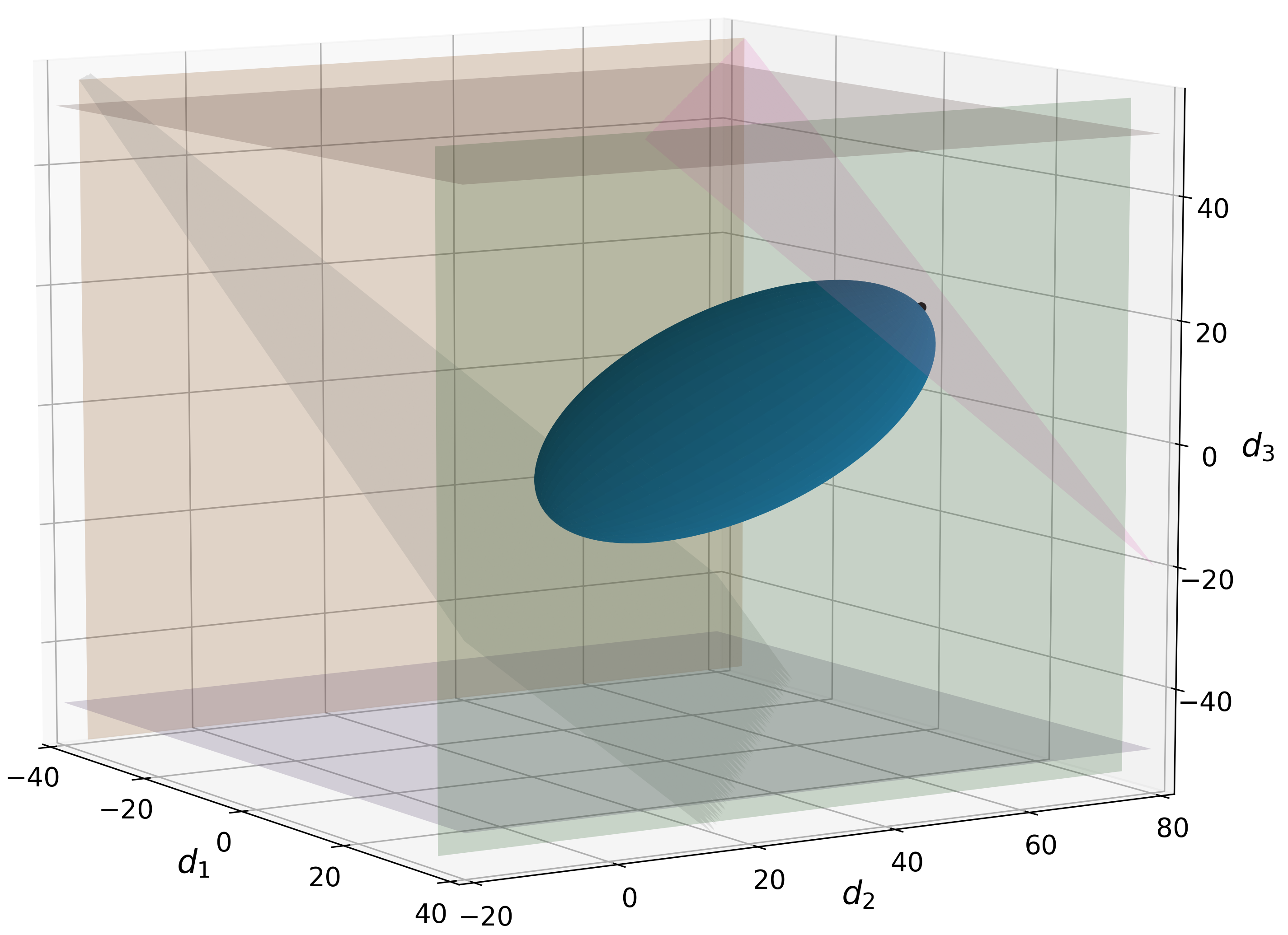}
		\caption{$\beta = 50$}
	\end{subfigure}
	\caption{Ellipsoidal sets with $\bar{\btheta} = \bar{\btheta}_{ac}$ and different covariance values $\beta$ (Design 1 of the three-node network).}
	\label{fig:3_node_covar_compare}
\end{figure}

We substituted the node balances into the capacity constraints to obtain a system of inequalities that are expressed solely in terms of $\btheta$. Thus, the uncertainty sets for Design 1 can be visualized in three dimensions. Figure \ref{fig:3_node_sets} depicts the uncertainty sets of Table \ref{tab:3node_results} with the center $\bar{\btheta}_{fc}$. Here, it is clear how intersecting the sets $T_{ellip}(\delta)$ and $T_{box}(\delta)$ with $\mathbb{R}^{n_{\theta}}_+$ creates more exotic regions and this affects the solution. Specifically, the intersected sets are scaled to a greater extent and they are limited by different constraints (in comparison to their non-intersected counterparts). This again reflects that the shape of $T(\delta)$ has a significant effect on the behavior of the flexibility index. Also, when Figure \ref{fig:ellip_fc} is juxtaposed with Figure \ref{fig:ellip_ac},  it is apparent how the analytic center is more neutrally positioned between the supplier capacity constraints (slanted planes) relative to the feasible center. This observation explains why the flexibility is identical for the intersected and non-intersected sets. 

We now show  $T_{ellip}(\delta)$ can be used to conservatively approximate $SF$ via the corresponding confidence level $\alpha^*$ and how correlation in $\btheta$ affects indexes $SF$ and $F$. Specifically, we consider the combination of instances for $\beta \in \{-40, 0, 50\}$ and $\bar{\btheta} \in \{\bar{\btheta}_{ac}, \bar{\btheta}_{fc}\}$. Table \ref{tab:3node_covar_compare} summarizes the solutions to these instances. Again, we observe that $F_{ellip}$ (and the corresponding confidence levels) mirror the trends observed with $SF$ at a significantly reduced computational cost. The correlation between the random variables clearly has an impact on $SF$, which  $F_{ellip}$ captures as well. On the other hand, index $F_{box}$ does not capture this behavior since it cannot account for parameter correlation. Also, it is apparent that the confidence levels associated with $\beta = -40$ provide a very tight lower bound for $SF$, whereas the one associated with $\beta = 50$ has a much larger gap. This behavior illustrates how the shape of the feasible region and that of the uncertainty set affect the tightness of the bound on $SF$ achieved by $\alpha^*$. 

\begin{table}[!htb]
	\caption{Results for Design 1 for various covariance values $\beta$.}
	\begin{center}
		\begin{tabular}{|c|ccccc|}
			\hline
			Center   & $\beta$ & $F_{ellip}$ & $\alpha^*$ (\%) & $SF$-MC (\%) & Solution Time (s) \\ \hline \hline
			         & -40     & 15.31       & 99.84           & 99.99        & 0.18 \\
			analytic & 0       & 8.93        & 96.97           & 99.71        & 0.21 \\
			         & 50      & 4.31  		 & 77.02           & 96.22        & 0.18 \\ \hline
			         & -40     & 15.31       & 99.84           & 99.99        & 0.14 \\
			feasible & 0       & 5.00        & 82.79           & 98.71        & 0.16 \\
			         & 50      & 2.41        & 50.85           & 93.49        & 0.15 \\ \hline
		\end{tabular}
	\end{center}
	\label{tab:3node_covar_compare}
\end{table}

Figure \ref{fig:3_node_covar_compare} depicts the uncertainty sets corresponding to the results of Table \ref{tab:3node_covar_compare}. We note how the ellipsoidal set with $\beta = -40$ is able to scale to a greater extent relative to the set associated with $\beta = 50$ because of its more favorable shape relative to the feasible region. This helps explain why $\alpha^*$ is able to provide the tightest lower bound for $SF$ with $\beta = -40$. 

\FloatBarrier

\subsubsection{Design Comparison}

To compare the three different designs described in Figure \ref{fig:3_node_designs}, we consider the results obtained with $\bar{\btheta} = \bar{\btheta}_{ac}$ (the analytic center for Design 1) and $\beta = 50$. Table \ref{tab:3node_design_compare} summarizes these results. The results indicate that Design 2 (which features decentralized suppliers) is less flexible than Design 1 (which utilizes a centralized supplier scheme). This is surprising, since intuition suggests that decentralization increases network flexibility. Our framework, however, shows that the converse is true in this particular instance. This behavior results from the correlation structure of the random variables. Furthermore, we observe that the arc added in Design 3 does not increase system flexibility (as one might also intuitively expect). As we discuss next, this behavior results from non-obvious interplays in the limiting constraints. This highlights how flexibility analysis can help guide system design.

\begin{table}[!htb]
	\caption{A summary of the flexibility index results using $F_{ellip}$ and $F_{box}$ to compare Designs 1, 2, and 3 of the three-node distribution network. }
	\begin{center}
		\begin{tabular}{|c|cccc|}
			\hline
			         & $F_{box}$ & $F_{ellip}$ & $\alpha^*$ (\%) & $SF$-MC (\%)  \\ \hline \hline
			Design 1 & 0.578     & 4.310       & 77.02           & 96.22 \\
			Design 2 & 0.578     & 3.612       & 69.35           & 94.32 \\
			Design 3 & 0.578     & 4.310       & 77.02           & 96.20 \\ \hline
		\end{tabular}
	\end{center}
	\label{tab:3node_design_compare}
\end{table}

\subsubsection{Limiting Constraint Ranking}

\begin{figure}[!htb]
	\centering
	\includegraphics[width=0.5\textwidth]{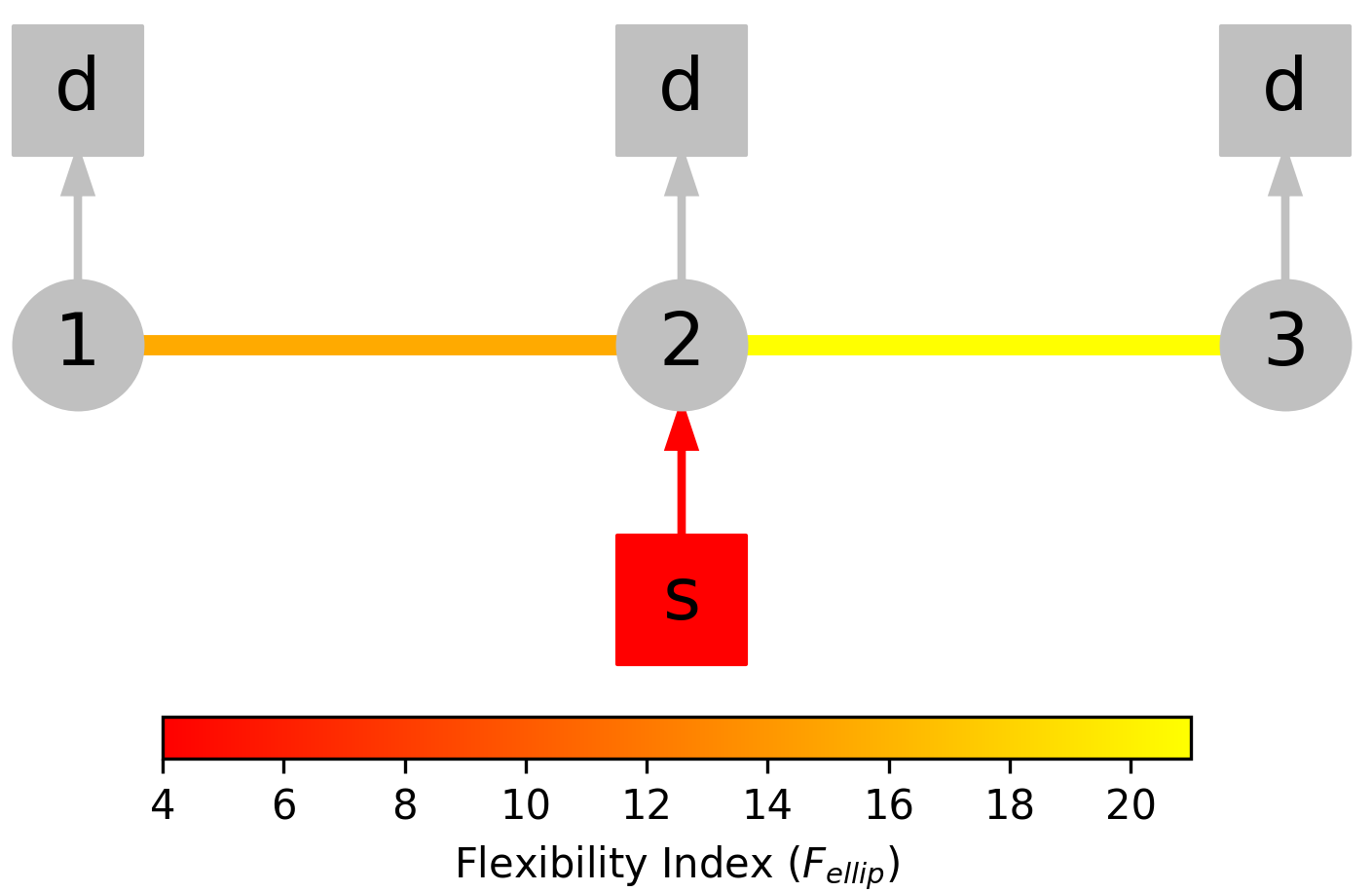}
	\caption{Limiting components for three-node network. Colors are proportional to the corresponding indexes $F_{ellip}$ (the smaller the value the most limiting the component is).}
	\label{fig:3_ranking}
\end{figure}

We use Design 1 to illustrate the limiting constraint ranking methodology described in Section \ref{sec:flex_rank}, letting $\bar{\btheta} = \bar{\btheta}_{ac}$ and $\beta = 50$ in conjunction with $T_{ellip}(\delta)$. Table \ref{tab:3node_ranking} summarizes the results. We observe that the supplier capacity limits system flexibility the most. This explains why adding an arc in Design 3 does not increase the flexibility of the network (i.e., the system is constrained by the inability to supply power and not by the inability to transport power). As a result, one can improve flexibility the most by increasing the supplier capacity (instead of adding an arc). Figure \ref{fig:3_ranking} depicts these results by shading the limiting design components in accordance with their corresponding index $F_{ellip}$. Here, components are colored according to the value of $F_{ellip}$ corresponding to their rank level (a lower value indicates that the component is more limiting). We thus see that supply is the most limiting component, the left arc is the second most limiting component, and the right arc is the third most limiting component.

\begin{table}[!htb]
	\caption{A summary of the constraint ranking results corresponding to Design 1 of the three-node distribution network using the set $T_{ellip}(\delta)$ with $\beta = 50$ and $\bar{\btheta} = \bar{\btheta}_{ac}$.}
	\begin{center}
		\begin{tabular}{|c|ccc|}
			\hline
			       & Active Constraint Multipliers        & $F_{ellip}$ & Flexibility Increase (\%) \\ \hline \hline
			Rank 1 & $\gamma^U_2$, $\gamma^L_2$           & 4.310       & $-$  \\
			Rank 2 & $\lambda^U_{2:1}$, $\lambda^L_{2:1}$ & 15.312      & 255.3 \\
			Rank 3 & $\lambda^U_{2:3}$, $\lambda^L_{2:3}$ & 20.833      & 383.4 \\ \hline
		\end{tabular}
	\end{center}
	\label{tab:3node_ranking}
\end{table}

\begin{figure}[!htb]
	\centering
	\begin{subfigure}[t]{.48\textwidth}
		\centering
		\includegraphics[width=\textwidth]{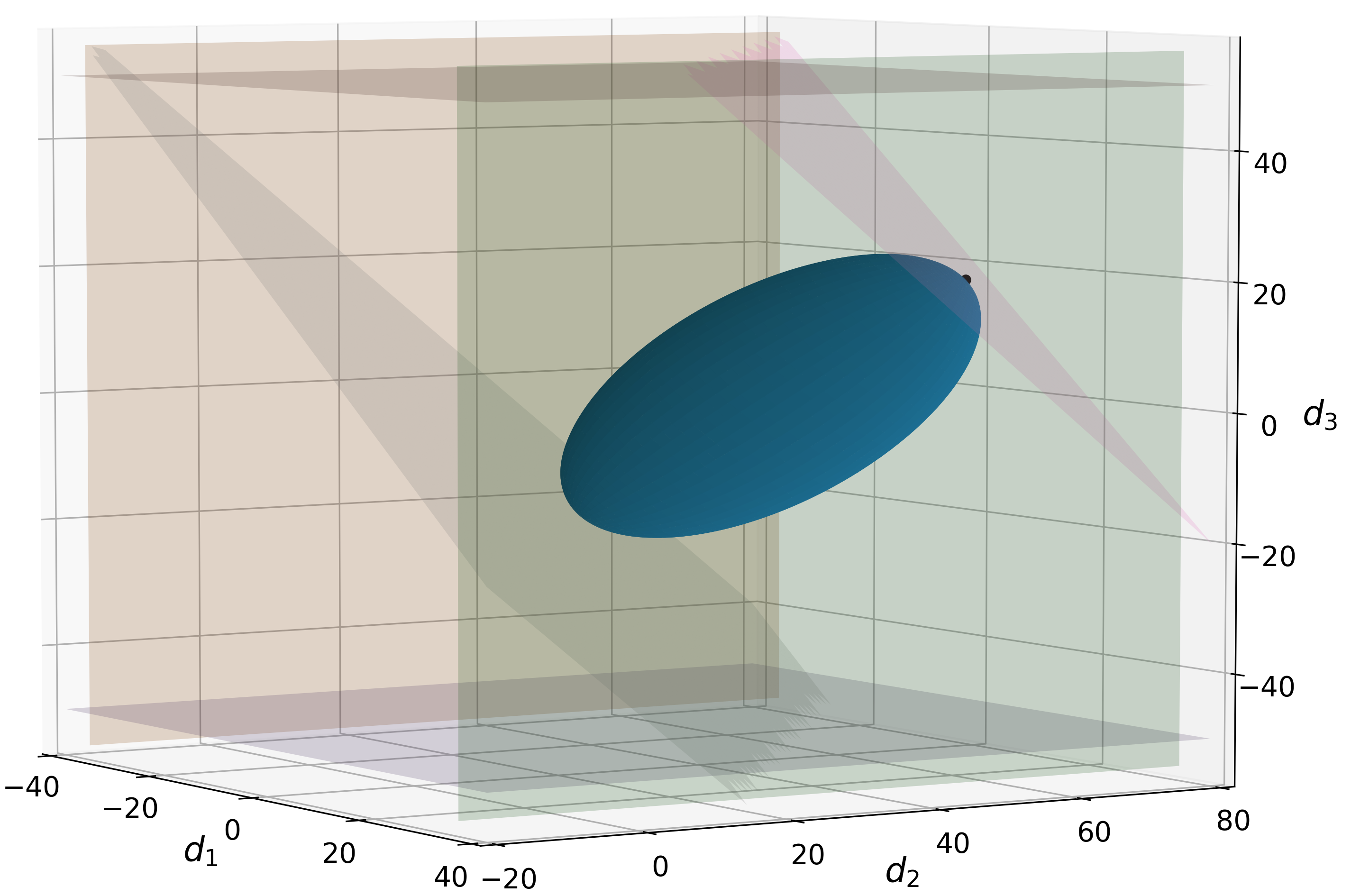}
		\caption{Rank 1}
	\end{subfigure}
	\quad
	\begin{subfigure}[t]{.48\textwidth}
		\centering
		\includegraphics[width=\textwidth]{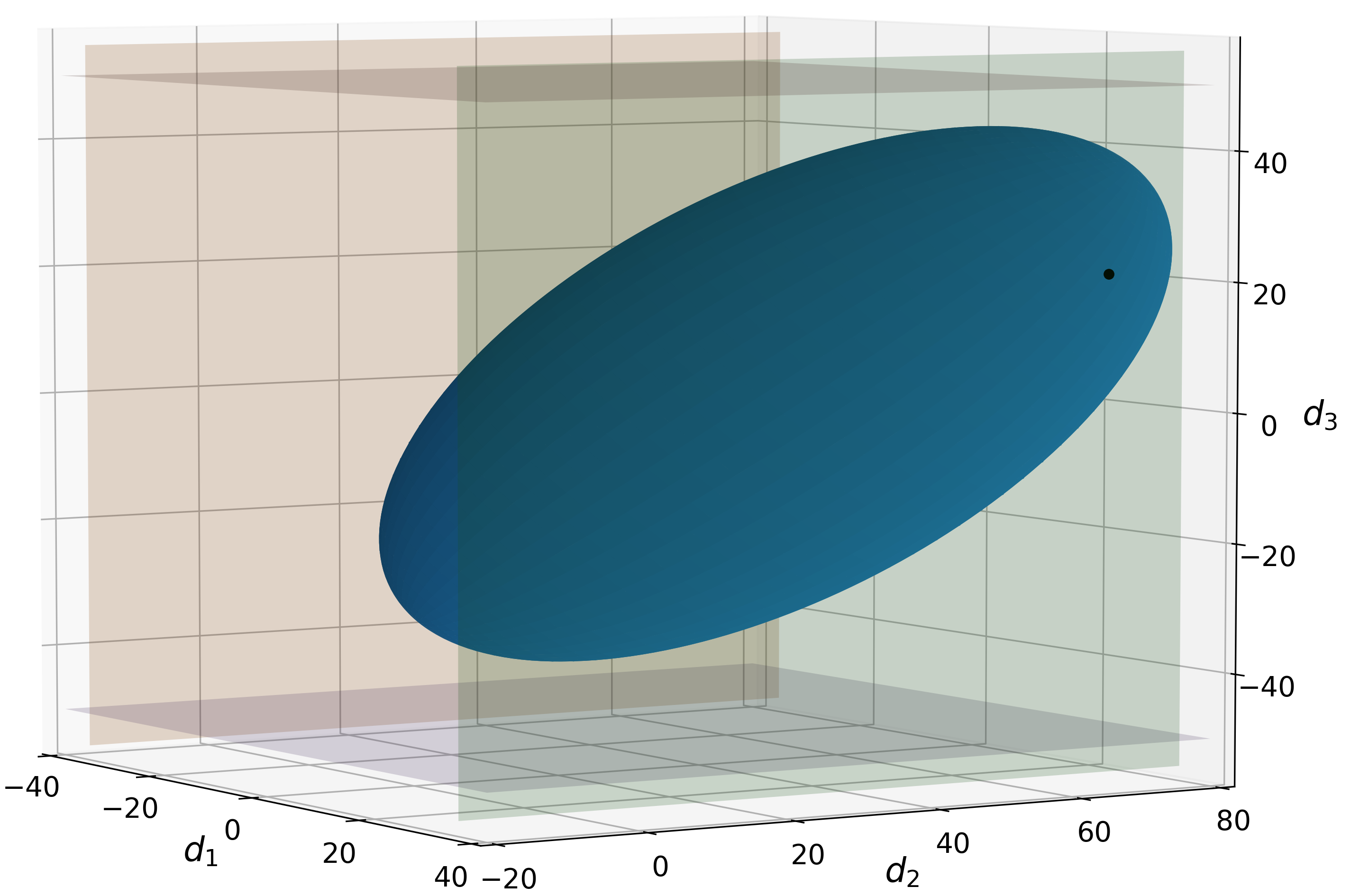}
		\caption{Rank 2}
	\end{subfigure} \\
	\vspace{1em}
	\begin{subfigure}[t]{.48\textwidth}
		\centering
		\includegraphics[width=\textwidth]{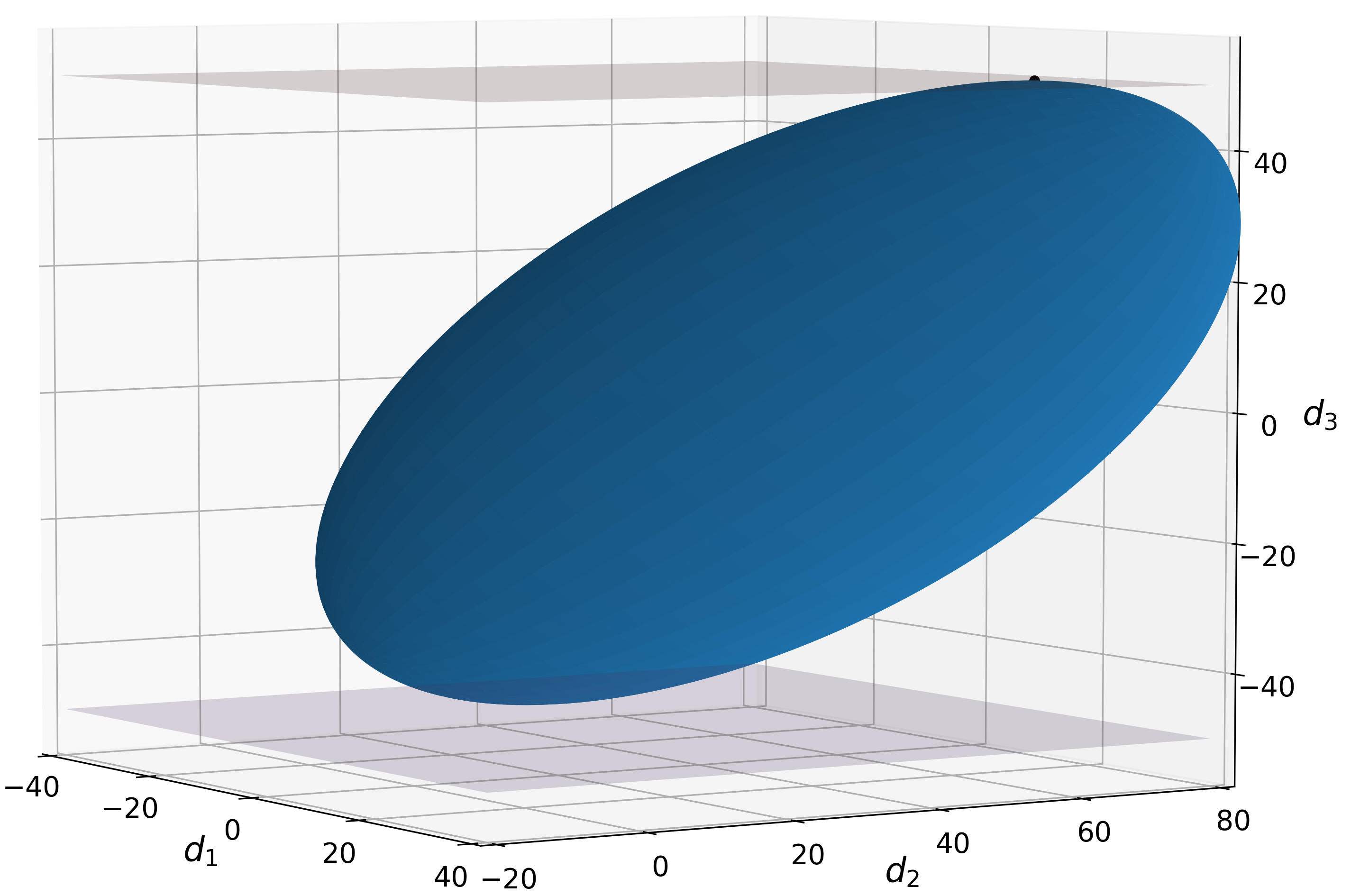}
		\caption{Rank 3}
	\end{subfigure}
	\caption{Limiting components for Design 1 of three-node network using $T_{ellip}(\delta)$ with $\bar{\btheta} = \bar{\btheta}_{ac}$ and $\beta = 50$.}
	\label{fig:3node_rank}
\end{figure}

Figure \ref{fig:3node_rank} depicts the solutions corresponding to the results of Table \ref{tab:3node_ranking}. We see that the size of $T_{ellip}(\delta)$ significantly increases as limiting constraints are removed. In this study, we found that two constraints simultaneously limit flexibility (corresponding to the maximum and minimum capacities of the limiting component).  This indicates that the uncertainty set touches the boundary of the feasible set at two locations.

\FloatBarrier

\subsection{IEEE 14-Node Network}

\begin{figure}[!htb]
	\centering
	\includegraphics[width=0.5\textwidth]{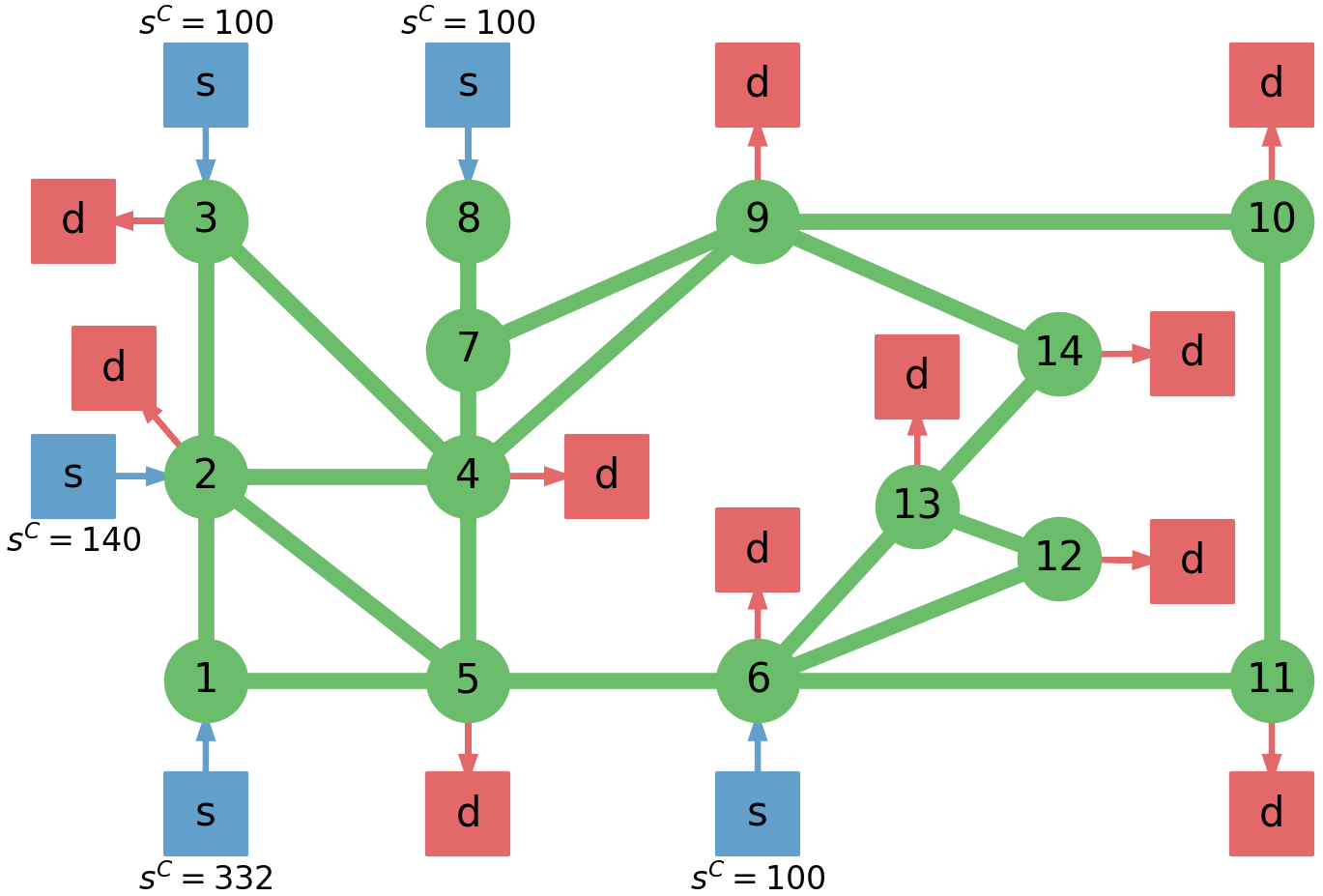}
	\caption{Schematic of the IEEE 14-node power system where the values $s^C$ are indicated and all the values $a^C = 100$.}
	\label{fig:14_diagram}
\end{figure}

We now consider the IEEE 14-node power network, originally provided in Dabbagchi in \cite{ieee14_origin}. The system data is obtained from MATPOWER. This test case does not provide arc capacities so we enforce $a^C = 100$ for all the arcs. This base design is labeled Design 1 and a schematic of this system is provided in Figure \ref{fig:14_diagram}. Design 2 is obtained by removing the arc that connects node 2 to node 5 and by removing the arc that connects node 9 to node 14. Design 3 uses a centralized power supply scheme, this is accomplished by removing all of the suppliers except for the two located on nodes 1 and 2, by increasing the capacity of each remaining supplier from 332 to 432 and from 140 to 340 respectively, and by increasing the capacity of each arc to 200. 

This system is subjected to a total of 10 uncertainty disturbances (the network demands). The demands are assumed to be $\btheta \sim \mathcal{N}(\bar{\btheta}, V_{\btheta})$, where $\bar{\btheta} = \bar{\btheta}_{fc} = (87.3,50.0,25.0,28.8,50.0,25.0,\allowbreak0,0,0,0,0)$ and $V_{\btheta} = 1200 \mathbb{I}$. We select the sets $T_{box+}(\delta)$ and $T_{ellip+}(\delta)$ (the demands are nonnegative). Each element of the hyperbox deviations $\Delta\btheta^-, \Delta\btheta^+$ is set to $3\sigma_i = 103.92$ (corresponding to $\bar{\btheta}_i \pm 3\sigma_i$ confidence bounds). Both of these sets are used to determine the flexibility index from Problem \eqref{eq:flex_mip} for each design. The stochastic flexibility index $SF$ is computed using 1,000,000 MC samples.

Table \ref{tab:ieee14_design_compare} summarizes these design comparison results. Index $SF$ shows that the flexibility of Design 1 is worsened with the removal of the two arcs, whereas the centralized supply modification improves flexibility. An observation is that $F_{ellip+}$ mirrors the behavior  $SF$ at a much lower computational cost. Specifically, $F_{ellip+}$ required less than 20 seconds to compute while computing $SF$ using MC sampling required 13,030 seconds (3.6 hours). Index $F_{box+}$ is not able to mirror the same behavior (it does not indicate that Design 2 is less flexible than Design 1 as would be expected). This analysis also shows that the higher capacity of the centralized supply scheme (Design 3) increases system flexibility relative to Design 1, which again is counter-intuitive. 

\begin{table}[!htb]
	\caption{Flexibility index results for various designs of IEEE-14 system.}
	\begin{center}
		\begin{tabular}{|c|ccc|}
			\hline
			         & $F_{box+}$ & $F_{ellip+}$ & $SF$-MC (\%)  \\ \hline \hline
			Design 1 & 0.327      & 10.594       & 97.84 \\
			Design 2 & 0.327      & 8.333        & 94.89 \\
			Design 3 & 0.404      & 12.244       & 99.91 \\ \hline
		\end{tabular}
	\end{center}
	\label{tab:ieee14_design_compare}
\end{table}

The limiting constraints of Design 1 are identified and ranked using $T_{ellip+}(\delta)$. The ranking results are shown in Table \ref{tab:ieee14_ranking}. The capacity constraints corresponding to four of the five suppliers along with three arc capacity constraints are identified as the most limiting ones. Again, this indicates that the uncertainty set touches the boundary of the feasible set at multiple points. We see that the next set of constraints only have an index that is 15.6\% larger meaning they are likely to also limit the system flexibility to a certain extent. Subsequent sets have significantly increased flexibility indexes (up to 214\% larger) and therefore have little effect on flexibility. 

\begin{table}[!htb]
	\caption{Ranking of limiting constraints for IEEE-14 system.}
	\begin{center}
		\begin{tabular}{|c|cccc|}
			\hline
		    & Active Constraints & $F_{ellip+}$ & Flexibility Increase (\%) & Solution Time (s)\\ \hline \hline
			Rank 1 & $\lambda_{1:2}^{U}$, $\lambda_{1:5}^{U}$, $\lambda_{6:11}^{U}$, $\gamma_{2}^{U}$, $\gamma_{3}^{U}$, $\gamma_{6}^{U}$, $\gamma_{8}^{U}$ & 10.594 & $-$ & 16.61 \\
			Rank 2 & $\lambda_{2:3}^{L}$, $\lambda_{2:5}^{L}$, $\lambda_{12:13}^{U}$, $\gamma_1^{L}$, $\gamma_2^{L}$, $\gamma_3^{L}$, $\gamma_6^{L}$, $\gamma_8^{L}$ & 12.243 & 15.6 & 3.76  \\
			Rank 3 & $\lambda_{2:3}^{U}$, $\lambda_{2:4}^{U}$, $\lambda_{3:4}^{U}$, $\lambda_{6:12}^{U}$, $\lambda_{6:13}^{U}$, $\lambda_{9:14}^{U}$ & 25.000 & 136.0 & 2.72  \\
			Rank 4 & $\lambda_{2:5}^{U}$, $\lambda_{5:6}^{L}$, $\lambda_{9:10}^{U}$, $\lambda_{9:14}^{L}$, $\lambda_{10:11}^{L}$, $\gamma_1^{U}$ & 33.333 & 214.6 & 0.14  \\ \hline          
		\end{tabular}
	\end{center}
	\label{tab:ieee14_ranking}
\end{table}

\begin{figure}[!htb]
	\centering
	\includegraphics[width=0.4\textwidth]{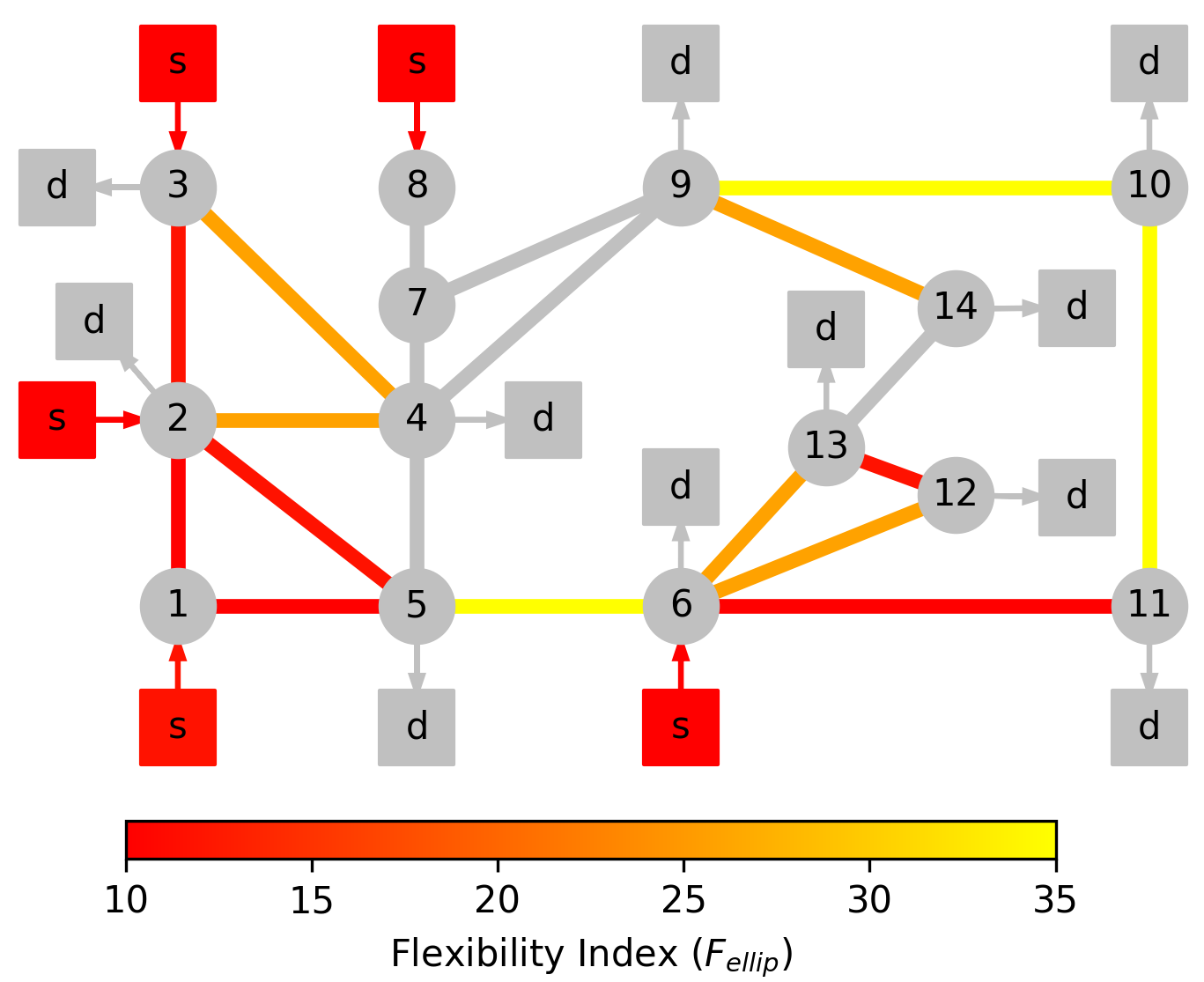}
	\caption{Limiting components for IEEE-14 power system. The colors are proportional to the corresponding indexes $F_{ellip+}$ (the smaller the value the most limiting the component is). }
	\label{fig:14_ranking}
\end{figure}

Figure \ref{fig:14_ranking} highlights the limiting components of Design 1. The limiting constraints corresponding to the active capacity constraints at each rank level are shaded according to the value of $F_{ellip+}$ (the lower the value, the most limiting the component is). We can see that the suppliers along with arcs attached to such suppliers limit flexibility the most.  The rest of the components follow non-intuitive patterns, reflecting complex dependencies due to the network topology. 

\FloatBarrier

\subsection{141-Node Network} \label{sec:141_node}

\begin{figure}[!htb]
	\centering
	\includegraphics[width=0.8\textwidth]{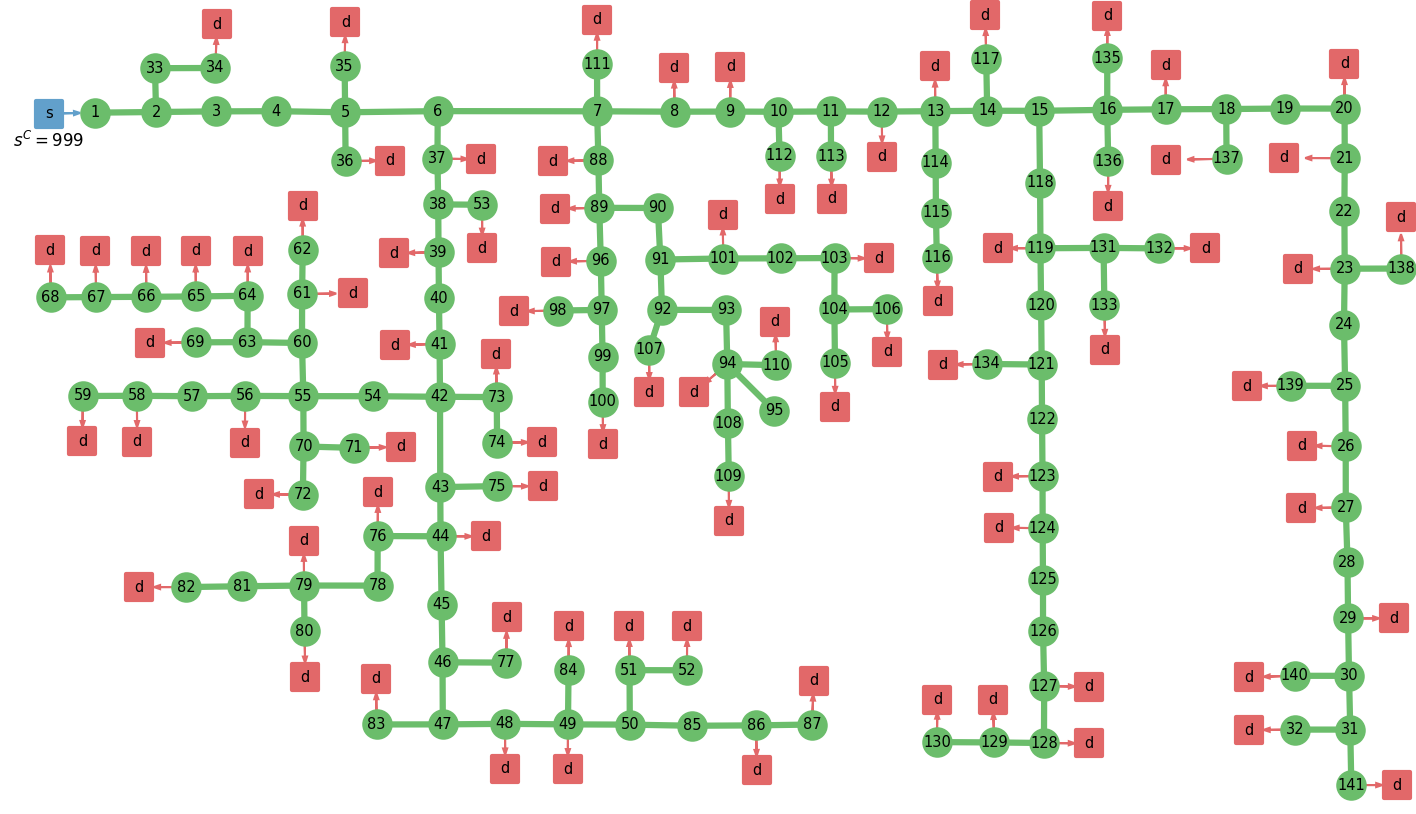}
	\caption{Schematic of the 141-node power distribution network.}
	\label{fig:141_diagram}
\end{figure}

Finally, we consider a 141-node power distribution network which corresponds to an urban area in Caracas, Venezuela and was originally developed in \cite{khodr2008maximum}. The network data is extracted from MATPOWER, but again no arc capacities are provided (we assume $a^C = 100$). Figure \ref{fig:141_diagram} provides a schematic of the network. 

This system is subjected to 84 uncertain disturbances (corresponding to the demands). The demands are assumed to be $\btheta \sim \mathcal{N}(\bar{\btheta}, V_{\btheta})$, where $\bar{\btheta} = \bar{\btheta}_{fc}$ and $V_{\btheta} = 100 \mathbb{I}$.  We use $T_{ellip}(\delta)$ in combination with Problem \eqref{eq:flex_mip} to rank all the inequality constraints.  This problem was solved 270 times while iteratively removing the active constraints. We found that several solutions have the same value of $F_{ellip}$. Thus, the active constraints corresponding to solutions with equivalent $F_{ellip}$ values are combined and assigned to the same ranks. A total of 44 rank levels were identified, the ranking data corresponding to the first 30 levels is provided in Appendix \ref{a:sup_mat}.

Figure \ref{fig:141_ranking} shows the limiting components corresponding to the 44 rank levels. We observe that the supplier is the most limiting component along with the arcs that form the spanning tree of the network (the path that connects all nodes in the network). The arcs become less relevant as they branch out towards the network boundaries.  We thus see that the flexibility framework reveals topological limitations of the network. We also highlight that the flexibility index can be computed for this large system in 20 seconds. In contrast, given the large number of random parameters, computing the $SF$ index using Monte Carlo sampling is impractical.

\begin{figure}[!htb]
	\centering
	\includegraphics[width=0.8\textwidth]{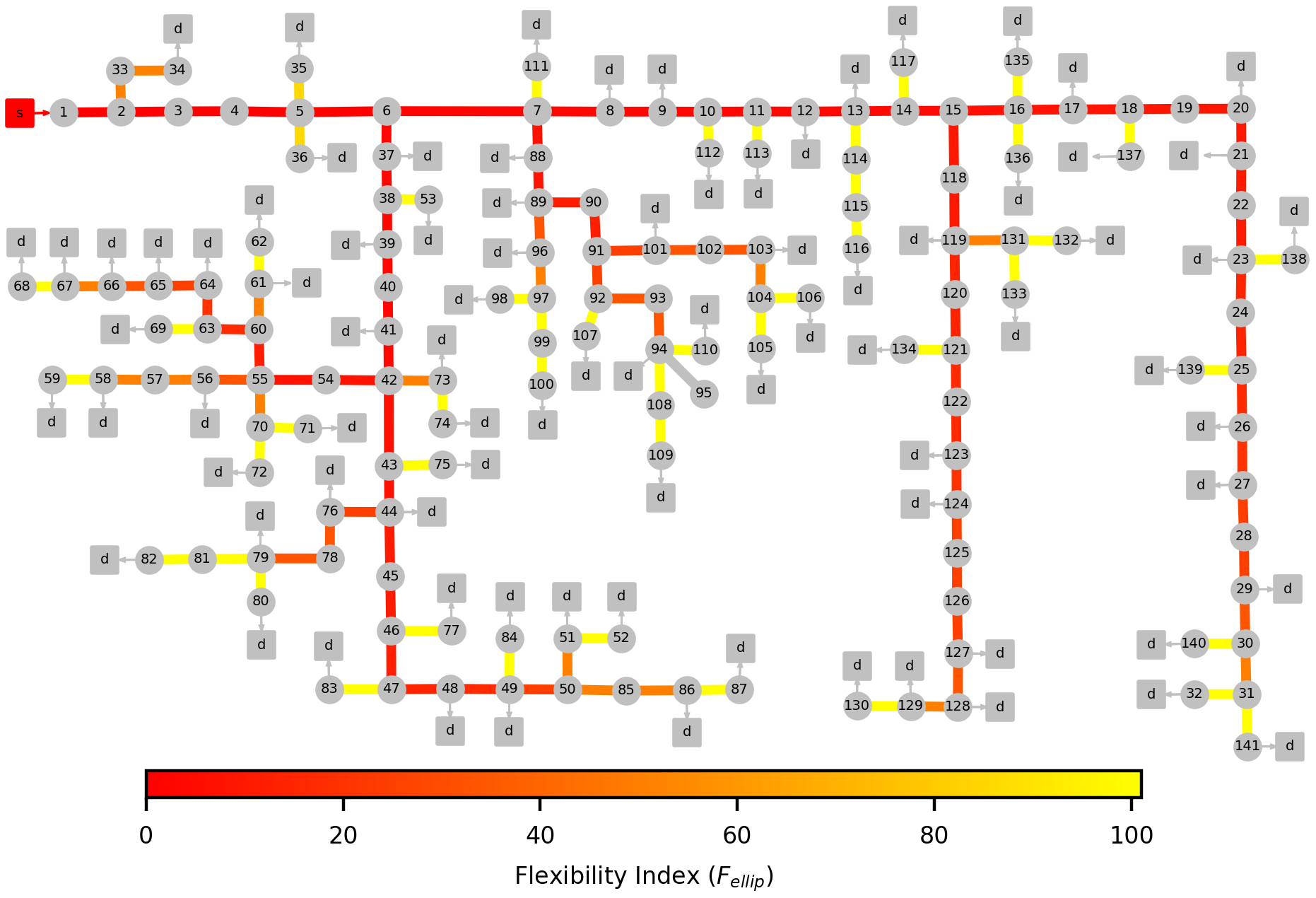}
	\caption{Limiting components for the 141-node power system. The colors are proportional to the corresponding indexes $F_{ellip}$ (the smaller the value the most limiting the component is).}
	\label{fig:141_ranking}
\end{figure}

\FloatBarrier

\section{Conclusions and Future Work}
We have presented a general framework to quantify and analyze system flexibility. This framework provides several new features that include: generalizing uncertainty sets to consider compositions of sets, finding a suitable nominal (center) point, and identifying and ranking limiting constraints. We have demonstrated that this framework is often able to emulate the behavior of the rigorous stochastic flexibility index problem at a significantly reduced computational cost. As part of future work, we are interested in developing techniques to handle nonlinear and dynamical systems such as reaction systems, heat exchanger networks, and transient power networks. We are also interested in developing techniques that do not require the solution of bilevel optimization problems. 

\section*{Acknowledgments}
This work was supported by the U.S. Department of Energy under grant DE-SC0014114.

\appendix
\section{Supplementary Material} \label{a:sup_mat}

The results for all of the instances discussed in Section \ref{sec:3node} are summarized in Tables \ref{tab:3node_norm_data} and \ref{tab:3node_pos_data}. Table \ref{tab:3node_norm_data} features the use of sets $T_{box}(\delta)$ and $T_{ellip}(\delta)$ for all designs, covariances, and centers considered in the three-node network example. Table \ref{tab:3node_pos_data} features the use of the sets $T_{box+}(\delta)$ and $T_{ellip+}(\delta)$ for the same system. Table \ref{tab:141_ranking} provides the first 30 rank levels corresponding to the results discussed in Section \ref{sec:141_node} for the 141-node network. 

\begin{table}[ht]
	\caption{Results for three-node distribution network using sets $T_{box}(\delta)$ and $T_{ellip}(\delta)$.}
	\begin{center}
		\begin{tabular}{|c|c|ccccccc|}
			\hline
			Design & Center   & $\beta$ & $F_{box}$ & Time (s) & $F_{ellip}$ & Time (s) & $\alpha^*$ (\%) & $SF$-MC (\%) \\ \hline \hline
			       &          & -40     & 0.578     & 0.03     & 15.31       & 0.18     & 99.84           & 99.99 \\
			       & analytic & 0       & 0.578     & 0.03     & 8.93        & 0.21     & 96.97           & 99.71 \\
			1      &          & 50      & 0.578     & 0.03     & 4.31  		 & 0.18     & 77.02           & 96.22 \\ \cline{2-9}
			       &          & -40     & 0.432     & 0.03     & 15.31       & 0.14     & 99.84           & 99.99 \\
			       & feasible & 0       & 0.432     & 0.03     & 5.00        & 0.16     & 82.79           & 98.71 \\
			       &          & 50      & 0.432     & 0.03     & 2.41        & 0.15     & 50.85           & 93.49 \\ \hline
		           &          & -40     & 0.578     & 0.08     & 3.61        & 0.21     & 69.35           & 96.82 \\
		           & analytic & 0       & 0.578     & 0.08     & 3.61        & 0.23     & 69.35           & 96.56 \\
		    2      &          & 50      & 0.578     & 0.08     & 3.61  		 & 0.19     & 69.35           & 94.32 \\ \cline{2-9}
		           &          & -40     & 0.432     & 0.03     & 3.61        & 0.15     & 69.35           & 96.82 \\
		           & feasible & 0       & 0.432     & 0.03     & 3.61        & 0.21     & 69.35           & 96.00 \\
		           &          & 50      & 0.432     & 0.03     & 2.41        & 0.20     & 50.85           & 92.67 \\ \hline
		           &          & -40     & 0.578     & 0.04     & 37.81       & 0.20     & 100.00          & 100.00 \\
		           & analytic & 0       & 0.578     & 0.04     & 8.93        & 0.16     & 96.97           & 99.72 \\
		    3      &          & 50      & 0.578     & 0.04     & 4.31  		 & 0.16     & 77.02           & 96.20 \\ \cline{2-9}
		           &          & -40     & 0.432     & 0.03     & 34.97       & 0.17     & 100.00          & 100.00 \\
		           & feasible & 0       & 0.432     & 0.03     & 5.00        & 0.16     & 82.79           & 98.72 \\
		           &          & 50      & 0.432     & 0.03     & 2.41        & 0.17     & 50.85           & 93.51 \\ \hline
		\end{tabular}
	\end{center}
	\label{tab:3node_norm_data}
\end{table}

\begin{table}[ht]
	\caption{Results for three-node distribution network using sets $T_{box+}(\delta)$ and $T_{ellip+} (\delta)$.}
	\begin{center}
		\begin{tabular}{|c|c|cccccc|}
			\hline
			Design & Center   & $\beta$ & $F_{box+}$ & Time (s) & $F_{ellip+}$ & Time (s) & $SF$-MC (\%) \\ \hline \hline
			       &          & -40     & 0.578     & 0.04     & 18.38       & 0.14       & 100.00 \\
			       & analytic & 0       & 0.578     & 0.04     & 8.93        & 0.10       & 99.44 \\
			1      &          & 50      & 0.578     & 0.04     & 4.31  		 & 0.13       & 94.33 \\ \cline{2-8}
			       &          & -40     & 0.723     & 0.03     & 18.38       & 0.11       & 100.00 \\
			       & feasible & 0       & 0.723     & 0.03     & 14.00       & 0.09       & 99.95 \\
			       &          & 50      & 0.273     & 0.03     & 6.76        & 0.09       & 98.60 \\ \hline
		           &          & -40     & 0.578     & 0.04     & 26.46       & 0.18       & 100.00 \\
			       & analytic & 0       & 0.578     & 0.04     & 8.82        & 0.16       & 99.34 \\
			2      &          & 50      & 0.578     & 0.04     & 4.31  		 & 0.13       & 94.23 \\ \cline{2-8}
			       &          & -40     & 0.723     & 0.03     & 26.46       & 0.12       & 100.00 \\
			       & feasible & 0       & 0.723     & 0.03     & 14.00       & 0.15       & 99.96 \\
			       &          & 50      & 0.723     & 0.03     & 6.76        & 0.12       & 98.60 \\ \hline
			       &          & -40     & 0.578     & 0.03     & 45.38       & 0.20       & 100.00 \\
			       & analytic & 0       & 0.578     & 0.03     & 8.93        & 0.11       & 99.45 \\
			3      &          & 50      & 0.578     & 0.03     & 4.31  		 & 0.06       & 94.34 \\ \cline{2-8}
			       &          & -40     & 0.723     & 0.03     & 47.19       & 0.09       & 100.00 \\
			       & feasible & 0       & 0.723     & 0.03     & 14.00       & 0.12       & 99.96 \\
			       &          & 50      & 0.723     & 0.03     & 6.76        & 0.09       & 98.60 \\ \hline
		\end{tabular}
	\end{center}
	\label{tab:3node_pos_data}
\end{table}

\FloatBarrier

\begin{table}[ht]
	\caption{Limiting constraints for 141-node network (determined using $T_{ellip}(\delta)$).}
	\begin{center}
		\begin{tabular}{|c|cccc|}
			\hline
			       & Active Constraint Multipliers & $F_{ellip}$ & Flexibility Increase (\%) & Solution Time (s)\\ \hline \hline
			Rank 1 & $\lambda_{1:2}^U$, $\gamma_{1}^L$ & 0.298 &  $-$  & 0.19 \\
			Rank 2 & $\lambda_{2:3}^U$, $\lambda_{4:5}^U$, $\lambda_{3:4}^U$ & 0.705  & 151.9 & 2.43 \\
			Rank 3 & $\lambda_{5:6}^U$  & 1.110  & 273.0   & 4.35 \\
			Rank 4 & $\lambda_{5:6}^L$   & 1.365 & 358.8 & 4.98 \\
			Rank 5 & $\lambda_{4:5}^L$, $\lambda_{2:3}^L$, $\lambda_{3:4}^L$ & 1.767 & 493.8 & 6.00 \\
			Rank 6 & $\lambda_{6:7}^U$   & 2.034 & 583.5 & 1.62 \\
			Rank 7 & $\lambda_{6:7}^L$   & 2.223 & 646.9 & 2.20 \\
			Rank 8 & $\lambda_{1:2}^L$    & 2.679 & 800.0   & 4.13 \\
			Rank 9 & $\lambda_{6:37}^U$   & 2.770  & 830.8 & 5.89 \\
			Rank 10 & $\lambda_{7:8}^U$   & 2.986 & 903.2 & 7.99 \\
			Rank 11 & $\lambda_{37:38}^{L,U}$ & 3.030  & 918.1 & 11.45 \\
			Rank 12 & $\lambda_{7:8}^L$   & 3.075 & 933.2 & 11.82 \\
			Rank 13 & $\lambda_{6:37}^L$   & 3.117 & 947.3 & 13.04 \\
			Rank 14 & $\lambda_{38:39}^{L,U}$, $\lambda_{8:9}^{L,U}$ & 3.125 & 949.9 & 25.67 \\
			Rank 15 & $\lambda_{40:41}^{L,U}$, $\lambda_{39:40}^{L,U}$, $\lambda_{9:10}^{L,U}$ & 3.226 & 983.9 & 44.29 \\
			Rank 16 & $\lambda_{41:42}^{L,U}$, $\lambda_{10:11}^{L,U}$ & 3.333 & 1020  & 18.95 \\
			Rank 17 & $\lambda_{11:12}^{L,U}$ & 3.448 & 1058.6 & 3.86 \\
			Rank 18 & $\lambda_{12:13}^{L,U}$ & 3.571 & 1099.9 & 4.64 \\
			Rank 19 & $\lambda_{13:14}^{L,U}$ & 3.846 & 1192.3 & 4.95 \\
			Rank 20 & $\lambda_{14:15}^{L,U}$ & 4.000     & 1243.9 & 5.08 \\
			Rank 21 & $\lambda_{15:16}^{L,U}$, $\lambda_{42:43}^{L,U}$ & 6.666 & 2139.9 & 14.41 \\
			Rank 22 & $\lambda_{43:44}^{L,U}$ & 7.143 & 2299.9 & 6.78 \\
			Rank 23 & $\lambda_{7:88}^{U}$  & 7.579 & 2446.5 & 6.10 \\
			Rank 24 & $\lambda_{16:17}^{L,U}$, $\lambda_{42:54}^{L,U}$, $\lambda_{54:55}^{L,U}$ & 7.692 & 2484.5 & 11.51 \\
			Rank 25 & $\lambda_{7:88}^{L}$  & 7.806 & 2522.9 & 4.06 \\
			Rank 26 & $\lambda_{17:18}^{L,U}$, $\lambda_{88:89}^{L,U}$ & 8.333 & 2699.9 & 9.83 \\
			Rank 27 & $\lambda_{19:20}^{L,U}$, $\lambda_{18:19}^{L,U}$ & 9.091 & 2954.5 & 3.78 \\
			Rank 28 & $\lambda_{20:21}^{L,U}$, $\lambda_{15:118}^{L,U}$, $\lambda_{118:119}^{L,U}$ & 10.000    & 3259.9 & 8.52 \\
			Rank 29 & $\lambda_{21:22}^{L,U}$, $\lambda_{22:23}^{L,U}$, $\lambda_{44:45}^{L,U}$, $\lambda_{45:46}^{L,U}$ & 11.111 & 3633.3 & 10.65 \\
			Rank 30 & $\lambda_{46:47}^{L,U}$, $\lambda_{55:60}^{L,U}$, $\lambda_{89:90}^{L,U}$, $\lambda_{90:91}^{L,U}$ & 12.500  & 4099.9 & 16.38 \\
			\hline      
		\end{tabular}
	\end{center}
	\label{tab:141_ranking}
\end{table}

\FloatBarrier

\bibliography{references}

\end{document}